\documentclass[12pt]{article}

\usepackage{amsfonts}
\usepackage{amsmath,amssymb,amsfonts}
\usepackage{color}
\usepackage{hyperref}
\def\squarebox#1{\hbox to #1{\hfill\vbox to #1{\vfill}}}
\newcommand{\qed}{\hspace*{\fill}
\vbox{\hrule\hbox{\vrule\squarebox{.667em}\vrule}\hrule}\smallskip}
\newtheorem{teorema}{Theorem}[section]
\newtheorem{lema}[teorema]{Lemma}
\newtheorem{corolario}[teorema]{Corollary}
\newtheorem{proposicao}[teorema]{Proposition}

\newenvironment{profe}{\noindent {\bf Proof:}}{\hfill $\qed $\vspace{10pt}}
\begin{document}

\title{Differentiability of Lyaupnov Exponents}
\author{Thiago F. Ferraiol\thanks{%
Address: Universidade Estadual de Maring\'{a}, Departamento de Matem\'{a}tica. Av. Colombo, 5790, Zona 7. 87020-900
Maring\'{a}, Paran\'{a}, Brasil. e-mail: tfferraiol@uem.br} \and Luiz A. B. San Martin\thanks{%
Supported by CNPq grant n$^{\mathrm{o}}$ 304982/2013-0 and FAPESP grant n$^{%
\mathrm{o}}$ 2012/17946-1.}\thanks{%
Address: Imecc - Unicamp, Departamento de Matem\'{a}tica. Rua S\'{e}rgio
Buarque de Holanda, 651, Cidade Universit\'{a}ria Zeferino Vaz. 13083-859
Campinas S\~{a}o Paulo, Brasil. e-mail: smartin@ime.unicamp.br}}
\maketitle

\begin{abstract}
We prove differentiability of certain linear combinations of the Lyapunov
spectra of a flow on a principal bundle of a semi-simple Lie group. The
specific linear combinations that yield differentiability are determined by
the finest Morse decomposition on the flag bundles. Differentiability is
taken with respect to a differentiable structure on the gauge group, which
is a Banach-Lie group.
\end{abstract}

\noindent%


\section{Introduction}

In this paper we prove differentiability of certain linear combinations of
the Lyapunov spectra of a continuous flow on a continuous principal bundle.
Our result is inspired on the theorem of Ruelle \cite{ruanl} that proves
differentiability of the top Lyapunov exponent for linear cocycles on vector
bundles under the assumption that the cocycle preserves a cone field.

Here we take cocycles with values on noncompact semi-simple Lie groups and
use the structure of the Morse decompositions on the respective flag bundles
to find linear combinations of the Lyapunov spectra that are differentiable
under perturbations of the generators of the cocycles. (Actually we use the
framework of principal bundles and perturb the flows by elements of the
gauge group which is a Banach-Lie group.)

As to the set up start with a continuous principal bundle $\pi :Q\rightarrow
X$ with compact Hausdorff base space $X$ and noncompact semi-simple Lie
group $G$ that acts on $Q$ on the right. An automorphism $\phi $ of $Q$ is a
homeomorphism $\phi :Q\rightarrow Q$ which is right invariant $\phi \left(
q\cdot g\right) =\phi \left( q\right) \cdot g$, $q\in Q$, $g\in G$, and
induces a homeomorphism $\tau :X\rightarrow X$ with $\pi \circ \phi =\tau
\circ \pi $. We consider the flow $\phi ^{n}$, $n\in \mathbb{N}$, generated
by $\phi $ and the corresponding flow $\tau ^{n}$ on $X$. We assume that
there is a probability measure $\nu $ on $X$, which is invariant by $\tau $
and has support $\mathrm{supp}\nu =X$.

The group $G$ admits an Iwasawa decomposition $G=KAN$ and a polar
decomposition $G=K\left( \mathrm{cl}A^{+}\right) K$ of $G$, with $%
A^{+}\subset A$, that are carried over to $Q$ (see \cite{alvsm} and \cite%
{smsec} for notation and details). In both decompositions the $A$-components
measure the exponential growth ratio of the iterations $\phi ^{n}\left(
q\right) $, $q\in Q$, for an automorphism $\phi $ of $Q$.

The logarithm of the $A$-component in the Iwasawa decomposition can be
factored to yield an additive cocycle $\mathsf{a}^{\phi }(n,\xi )\in 
\mathfrak{a}$ over the flow on the associated maximal flag bundle $\mathbb{E}%
=Q\times _{G}\mathbb{F}$ with values in the Lie algebra $\mathfrak{a}$ of
the abelian group $A$. It gives rise to the $\mathfrak{a}$-Lyapunov
exponents 
\begin{equation*}
\lambda \left( \xi \right) =\lim \frac{1}{n}\mathsf{a}^{\phi }(n,\xi )\in 
\mathfrak{a}\qquad \xi \in \mathbb{E},
\end{equation*}%
whose almost surely existence are ensured by the Multiplicative Ergodic
Theorem (see \cite{alvsm} for an approach in this context of semi-simple
groups).

On the other hand the projection onto $\mathrm{cl}A^{+}$-component against
the polar decomposition yields ultimately a subadditive cocycle $\omega
\circ \mathsf{a}_{\phi }^{+}(n,x)\in \mathrm{cl}\mathfrak{a}^{+}\subset 
\mathfrak{a}$ over the flow on the base space $X$, after composing with a
dominant weight $\omega $. The subadditive ergodic theorem applyied to the
several weights $\omega $ gives rise to the polar exponent 
\begin{equation*}
H_{\phi }^{+}\left( x\right) =\lim \frac{1}{n}\mathsf{a}_{\phi }^{+}(n,x)\in 
\mathrm{cl}\mathfrak{a}^{+}
\end{equation*}%
defined for $\nu $-almost all $x\in X$.

The vector exponents $\lambda \left( \xi \right) $, $\xi \in \mathbb{E}$,
and $H_{\phi }^{+}\left( x\right) $, $x\in X$, are related by the fact that
if $\lambda \left( \xi \right) $ is defined for some $\xi $ in the fibre
over $x$ then it is defined for every $\xi $ in this fiber and the possible
values are $\lambda \left( \xi \right) =wH_{\phi }^{+}\left( x\right) $ with 
$w$ running through the Weyl group of $G$ (see \cite{alvsm}).

We call $H_{\phi }^{+}\left( x\right) $ the $\mathfrak{a}$-Lyapunov spectrum
of $\phi $ with respect to $\nu $.

Our purpose in this paper is \ to prove that for suitable $\omega \in 
\mathfrak{a}^{\ast }$ the mean Lyapunov spectra 
\begin{equation}
\gamma \phi \in \mathcal{G}\phi \mapsto \omega \int H_{\gamma \phi
}^{+}\left( x\right) \nu \left( dx\right) \in \mathbb{R}  \label{formapanlit}
\end{equation}%
depends smoothly of $\gamma $ (at $\gamma =\mathrm{id}$) when $\gamma $ runs
through the gauge group $\mathcal{G}=\mathcal{G}\left( Q\right) $ of $Q$.

The linear maps $\omega \in \mathfrak{a}^{\ast }$ that enter in (\ref%
{formapanlit}) are given by dynamical properties of $\phi $. To simplify
notation we describe them in a trivial bundle $Q=X\times G$. The flow $\phi
^{n}$ on $X\times G$ induces flows on the flag bundles $\mathbb{E}_{\Theta
}=X\times \mathbb{F}_{\Theta }$ where $\mathbb{F}_{\Theta }$ is a flag
manifold of $G$. The Morse decompositions in such flag bundles were studied
in \cite{smbflow} and \cite{msm} (see also Selgrade \cite{selgrade} for a
different approach specific to projective bundles and Colonius-Kliemann for
bundles of flags of subspaces of $\mathbb{R}^{d}$). In \cite{smsec} the
Morse decompositions are related to spectra.

By the results of \cite{smbflow} and \cite{msm} in any flag bundle there
exists a finest Morse decomposition. It is given algebraically as follows:
For $H\in \mathrm{cl}\mathfrak{a}^{+}$, write $h=\exp H$ and $\mathbb{A}%
\left( H\right) =\{ghg^{-1}:g\in G\}$. Then there exists $H_{\mathrm{Mo}%
}\left( \phi \right) \in \mathrm{cl}\mathfrak{a}^{+}$ and a continuous map $%
\chi _{\phi }:X\rightarrow \mathbb{A}\left( H_{\mathrm{Mo}}\left( \phi
\right) \right) $ such that in any flag bundle $X\times \mathbb{F}_{\Theta }$
the Morse decomposition is given fiberwise as 
\begin{equation*}
\bigcup_{x\in X}\{x\}\times \mathrm{fix}_{\Theta }\left( \chi _{\phi }\left(
x\right) \right)
\end{equation*}%
where $\mathrm{fix}_{\Theta }\left( \chi _{\phi }\left( x\right) \right) $
is the set of fixed points of $\chi _{\phi }\left( x\right) $ in the flag
manifold $\mathbb{F}_{\Theta }$ (see \cite{msm}, Theorem 7.5). The element $%
H_{\mathrm{Mo}}\left( \phi \right) \in \mathrm{cl}\mathfrak{a}^{+}$ itself
defines a flag manifold which depends only on the set%
\begin{equation*}
\Theta _{\mathrm{Mo}}\left( \phi \right) =\{\alpha \in \Sigma :\alpha \left(
H_{\mathrm{Mo}}\left( \phi \right) \right) =0\},
\end{equation*}%
where $\Sigma $ is the set of simple roots. We call $\mathbb{F}_{\Theta _{%
\mathrm{Mo}}\left( \phi \right) }$ (or just $\Theta _{\mathrm{Mo}}\left(
\phi \right) $) the flag type of $\phi $.

Then the linear maps $\omega \in \mathfrak{a}^{\ast }$ such that (\ref%
{formapanlit}) is smooth are obtained from the flag type $\Theta _{\mathrm{Mo%
}}\left( \phi \right) $ as follows: Write the set of simple roots as 
\begin{equation*}
\Sigma =\{\alpha _{1},\ldots ,\alpha _{l}\}.
\end{equation*}%
The corresponding set of fundamental weights $\Omega =\{\omega _{1},\ldots
,\omega _{l}\}$ is defined dually by 
\begin{equation*}
\frac{2\langle \omega _{i},\alpha _{j}\rangle }{\langle \alpha _{j},\alpha
_{j}\rangle }=\delta _{ij}.
\end{equation*}%
Let $\Omega _{\Theta _{\mathrm{Mo}}\left( \phi \right) }\subset \Omega $ be
the subset of those $\omega _{i}$ having the same index $i$ of a $\alpha
_{i}\in \Theta _{\mathrm{Mo}}\left( \phi \right) $.

Then our main result Theorem \ref{teoprinc} says that (\ref{formapanlit}) is
differentiable if $\omega $ is a linear combination of $\Omega \setminus
\Omega _{\Theta _{\mathrm{Mo}}\left( \phi \right) }$.

The example of $G=\mathrm{Sl}\left( d,\mathbb{R}\right) $ is enlightening.
In this case $H_{\mathrm{Mo}}\left( \phi \right) $ is a diagonal matrix%
\begin{equation*}
H_{\mathrm{Mo}}\left( \phi \right) =\mathrm{diag}\{a_{1}\geq a_{2}\geq
\cdots \geq a_{d}\}
\end{equation*}%
and the flag type $\Theta _{\mathrm{Mo}}\left( \phi \right) $ is read off
from the multiplicites of the eigenvalues of $H_{\mathrm{Mo}}\left( \phi
\right) $, since the set of simple roots is 
\begin{equation*}
\Sigma =\{\lambda _{1}-\lambda _{2},\ldots ,\lambda _{d-1}-\lambda _{d}\}
\end{equation*}%
where $\lambda _{i}\left( H\right) $ is the $i$-th entry of the diagonal
matrix $H$. The corresponding fundamental weights are $\omega _{i}=\lambda
_{1}+\cdots +\lambda _{i}$, $1\leq i\leq d-1$. It follows that (\ref%
{formapanlit}) is differentiable if $\omega $ is a linear combination of the
sums $\omega _{i}=\lambda _{1}+\cdots +\lambda _{i}$ such that $%
a_{i}>a_{i+1} $.

As in \cite{ruanl} our method of proof exploits the Implicit Function
Theorem for Banach manifolds. From the Morse decompositions we get a section 
$\sigma $ of the specific flag bundle $\mathbb{E}_{\Theta _{\mathrm{Mo}%
}\left( \phi \right) }$ defined by the flag type of $\phi $. This section is
the attractor Morse component of $\phi $ in $\mathbb{E}_{\Theta _{\mathrm{Mo}%
}\left( \phi \right) }$. The part of the spetrum given by $\omega \in 
\mathfrak{a}^{\ast }$ can be writen as the integral of $\omega \circ \mathsf{%
a}^{\phi }(1,\xi )$ with respect to the measure $\mu =\sigma _{\ast }\nu $
on $\mathbb{E}_{\Theta _{\mathrm{Mo}}\left( \phi \right) }$. With the aid of
the Implicit Function Theorem we can parametrize differentiably invariant
measures $\mu _{\gamma }=\left( \sigma _{\gamma }\right) _{\ast }\nu $ as
well as integrands $\omega \circ \mathsf{a}^{\gamma \phi }(1,\xi )$ for $%
\gamma \in \mathcal{G}$ close to the identity and hence get smoothness of
the spectra since they are integrals with respect to the invariant measures
on $\mathbb{E}_{\Theta _{\mathrm{Mo}}\left( \phi \right) }$.

As applications of our main result we consider in Section \ref{secsemi}
cocycles with values in semigroups whose flag types are known. It can be
proved that the flag type $\Theta _{\mathrm{Mo}}\left( \phi \right) $ of the
flow is contained in the flag type of the semigroup (see Proposition \ref%
{propflagcontido}), so that the latter furnishes some $\omega $ for which
differentiability holds. An example of such semigroup is the one considered
in \cite{ruanl}, namely the semigroup $S_{W}\subset \mathrm{Sl}\left( d,%
\mathbb{R}\right) $ leaving invariant a pointed and generating cone $%
W\subset \mathbb{R}^{d}$. Its flag type is known to be the projective space
hence we get differentiability for $\omega =\lambda _{1}$ as in \cite{ruanl}.

Other examples of semigroups are worked out in Section \ref{secsemi}. One of
them is the semigroup of totally positive matrices whose flag type is the
empty set which corresponds to the maximal flag manifold. In this extremal
case we get smoothness of the whole spectra.

We conclude this introduction with the following comments.

First we work throughout with semi-simple Lie groups in order to facilitate
the algebraic technicalities. However the results are easily extended to
reductive Lie groups (like $\mathrm{Gl}\left( d,\mathbb{R}\right) $) since
the central component of the cocycle $\mathsf{a}^{\gamma \phi }(1,\xi )$
descends to a cocycle on the base space (see \cite{alvsm}, Section 3.3)
whose integral with respect to $\nu $ is differentiable.

When the measure $\nu $ on the base space $X$ is ergodic then the
Multiplicative Ergodic Theorem provides a measurable section $\chi _{\mathrm{%
Ly}}$ on the bundle with fiber $\mathbb{A}\left( H_{\mathrm{Ly}}\left( \phi
\right) \right) $ of conjugates of $h_{\mathrm{Ly}}=\exp H_{\mathrm{Ly}%
}\left( \phi \right) $ for some $H_{\mathrm{Ly}}\left( \phi \right) \in 
\mathrm{cl}\mathfrak{a}^{+}$. The fixed point sets of $\chi _{\mathrm{Ly}%
}\left( x\right) $, $x\in X$, give the Oseledets type decomposition of the
flag bundles that describes the level sets of the Lyapunov exponents (see 
\cite{alvsm}). This $H_{\mathrm{Ly}}\left( \phi \right) $ yields a flag type 
$\Theta _{\mathrm{Ly}}\left( \phi \right) $ analogous to the flag type for
the Morse decomposition. Since the Oseledets decomposition is contained in
the Morse decomposition (see \cite{alvsm1}) we have $\Theta _{\mathrm{Ly}%
}\left( \phi \right) \subset \Theta _{\mathrm{Mo}}\left( \phi \right) $. In
general this inclusion is proper. As a consequence of our differentiability
result it was proved in \cite{alvsm1} that the whole spectra is continuous
if $\Theta _{\mathrm{Ly}}\left( \phi \right) =\Theta _{\mathrm{Mo}}\left(
\phi \right) $.

Our result proves differentiability of the Lyapunov spectra when the flow $%
\phi $ is perturbed by elements of the gauge group $\mathcal{G}$. However
the properties of $\mathcal{G}$ needed in the proof may be satisfied by
other groups of automorphisms of $Q$. In Section \ref{secgener} we
abstracted these properties of $\mathcal{G}$ into five conditions to be
satisfied by a group $\mathcal{H}$ so that the proof still works with $%
\mathcal{H}$ in place of $\mathcal{G}$. When these conditions are satisfied
the Lyapunov spectra behaves smoothly when perturbing $\phi $ by elements of 
$\mathcal{H}$.

\section{Additive cocycles over flag bundles}

We recall here the construction of the $\mathfrak{a}$-cocycle yielding the
vector valued Lyapunov exponents (see \cite{alvsm}, \cite{alvsm1}, \cite%
{smsec} for further details).

Let $Q\rightarrow X$ be a principal bundle whose structural group $G$ is a
reductive Lie group. The right action of $G$ on $Q$ is denoted by $q\mapsto
q\cdot g$, $q\in Q$, $g\in G$. We endow $G$ with an Iwasawa decomposition $%
G=KAN$, that is kept fixed throughout.

The bundle $Q$ has a reduction to a subbundle $R\subset Q$ which is a
principal bundle with structural group $K$ in such a way that $Q=R\times AN$%
. This means that any $q\in Q$ decomposes uniquely as 
\begin{equation*}
q=r\cdot hn,\quad r\in R,\quad hn\in AN.
\end{equation*}%
We let 
\begin{equation*}
\mathsf{R}:Q\rightarrow R,\quad q\mapsto r,\quad \qquad \mathsf{A}%
:Q\rightarrow A,\quad q\mapsto h,
\end{equation*}%
be the ensuing projections. They satisfy the following properties:

\begin{enumerate}
\item $\mathsf{R}(r)=r$, $\mathsf{A}(r)=1$, when $r\in R$,

\item $\mathsf{R}(q\cdot p)=\mathsf{R}(q)m$, $\mathsf{A}(q\cdot p)=\mathsf{A}%
(q)h$, when $q\in Q$, $p=mhn\in P=MA\!N$. In particular, $\mathsf{A}(r\cdot
p)=h$.
\end{enumerate}

In what follows we write for $q \in Q$, 
\begin{equation*}
\mathsf{a}\left( q\right) =\log \mathsf{A}\left( q\right) \in \mathfrak{a}.
\end{equation*}

We use the same notation for the Iwasawa decomposition of $g\in G$, namely, $%
\mathsf{a}\left( g\right) =\log h$ if $g=uhn\in KAN$. For future reference
we note that by the second of the above properties we have 
\begin{equation}
\mathsf{a}\left( q\cdot p\right) =\mathsf{a}\left( q\right) +\mathsf{a}%
\left( p\right) ,\quad p\in P.  \label{forcocaditivo}
\end{equation}

Now let $\phi $ be an automorphism of $Q$. It induces the map 
\begin{equation*}
\phi ^{R}:r\in R\mapsto \mathsf{R}(\phi (r))\in R,
\end{equation*}%
which satisfies $\phi ^{R}\circ \psi ^{R}=(\varphi \circ \psi )^{R}$ if $%
\psi $ is another automorphism. Hence the flow defined by $\phi $ induces
the flow $\phi ^{R}$ on $R$. The map 
\begin{equation*}
\mathsf{a}^{\phi }(n,r)=\mathsf{a}(\phi ^{n}(r))
\end{equation*}%
is an additive cocycle over $\phi ^{R}$, that is, 
\begin{equation*}
\mathsf{a}^{\phi }(n+m,r)=\mathsf{a}^{\phi }(n,\left( \phi ^{R}\right)
^{m}(r))+\mathsf{a}^{\phi }(m,r).
\end{equation*}

This cocycle over $R$ factors to a cocycle over the flag bundle $\mathbb{E}%
=Q\times _{G}\mathbb{F}$, denoted by the same symbol (see \cite{alvsm},
Section 3). \ The $\mathfrak{a}$-Lyapunov exponent of $\phi $ in the
direction of $\xi \in \mathbb{E}$ is defined by 
\begin{equation*}
\lambda \left( \xi \right) =\lim_{t\rightarrow +\infty }\frac{1}{n}\mathsf{a}%
^{\phi }\left( n,\xi \right) \in \mathfrak{a},
\end{equation*}%
if the limit exists.

The cocycle $\mathsf{a}^{\phi }(n,\xi )$ is defined only on the maximal flag
bundle $\mathbb{E}$. It does not factor to a cocycle on a smaller flag
bundle $\mathbb{E}_{\Theta }$. However by composing it with suitable linear
maps $\omega \in \mathfrak{a}^{\ast }$ we get cocycles $\mathsf{a}_{\omega
}^{\phi }\left( n,\xi \right) =\omega \circ \mathsf{a}^{\phi }\left( n,\xi
\right) $ that can be factored to a cocycle in $\mathbb{E}_{\Theta }$. As
checked in \cite{alvsm}, Section 6, if $\omega $ belongs to the span of $%
\Omega \setminus \Omega _{\Theta }$ (see the introduction) then 
\begin{equation*}
\mathsf{a}_{\omega }^{\phi }\left( n,r\right) =\omega \circ \mathsf{a}^{\phi
}\left( n,r\right)
\end{equation*}%
is constant along the fibers of $\mathbb{E}\rightarrow \mathbb{E}_{\Theta }$
and hence defines a cocycle over the flow induced by $\phi $ on $\mathbb{E}%
_{\Theta }$. Note that by definition $\mathrm{span}\left( \Omega \setminus
\Omega _{\Theta }\right) $ is the annihilator $\mathfrak{a}\left( \Theta
\right) ^{\circ }$ of \ref{foradeteta} 
\begin{equation}
\mathfrak{a}\left( \Theta \right) =\{H\in \mathfrak{a}:\alpha \left(
H\right) =0,~\alpha \in \Theta \}.  \label{foradeteta}
\end{equation}

\section{ Manifold structures}

In this section we review some results about differentiable structures on
the gauge group and on the space of sections of associated bundles (see \cite%
{ferr}, Palais \cite{palais}, Wockel \cite{woc}).

Let $\pi :Q\rightarrow X$ be a principal bundle with group $G$ and compact
base space $X$. We assume that there exists a reduction $\pi :R\rightarrow X$
with compact group $K$.

If $F$ is a manifold acted differentiably on the left by $G$ we write $\pi
:Q\times _{G}F\rightarrow X$ for the associated bundle whose elements are
equivalence classes $q\cdot \xi $ of the equivalence relation on $Q\times F$%
, $\left( q,\xi \right) \sim \left( qg,g^{-1}\xi \right) $, $g\in G$. We
denote its space of continuous sections by $\Gamma \left( Q\times
_{G}F\right) $. For simplicity we assume throughout that $G$ acts
transitively on $F$, that is, $F=G/L$ for a closed subgroup $L\subset G$. We
assume that $K$ acts transitively on $F$, so that $F=G/L=K/\left( L\cap
K\right) $.

We topologize $\Gamma \left( Q\times _{a}F\right) $ in two different and
equivalent ways (we leave implicit that $\Gamma \left( Q\times _{a}F\right)
\neq \emptyset $). In both cases $E=Q\times _{G}F$ is viewed as the
associated bundle $R\times _{K}F\approx Q\times _{G}F$.

First let $\langle \cdot ,\cdot \rangle $ be a $K$-invariant Riemannian
metric in $F$ and denote by $d_{0}$ the corresponding distance. Since $d_{0}$
is $K$-invariant, $d\left( p\cdot z,p\cdot w\right) =d_{0}\left( z,w\right) $
defines a distance in the fiber of $R\times _{K}F$ over $\pi \left( p\right) 
$. Then 
\begin{equation*}
\rho \left( \sigma ,\tau \right) =\sup_{x\in X}d\left( \sigma \left(
x\right) ,\tau \left( x\right) \right) \qquad \sigma ,\tau \in \Gamma \left(
R\times _{K}F\right)
\end{equation*}%
is a distance in $\Gamma \left( R\times _{K}F\right) $.

Alternatively $\Gamma \left( Q\times _{G}F\right) $ is in bijection with the
space $\mathcal{C}_{\mathrm{eq}}\left( R,F\right) $ of continuous
equivariant maps $f:R\rightarrow F$, $f\left( qg\right) =g\cdot f\left(
q\right) $. The bijection associates $\sigma \in \Gamma \left( Q\times
_{G}F\right) $ to $f_{\sigma }\in \mathcal{C}_{\mathrm{eq}}\left( R,F\right) 
$ defined by $\sigma \left( x\right) =q\cdot f_{\sigma }\left( q\right) $,
for any $q$ with $\pi \left( q\right) =x$. In $\mathcal{C}_{\mathrm{eq}%
}\left( R,F\right) $ there is a distance defined by $\rho \left( f,g\right)
=\sup_{q\in R}d_{0}\left( f\left( q\right) ,g\left( q\right) \right) $. The
map $\sigma \mapsto f_{\sigma }$ is an isometry.

The manifold structure on $\Gamma E$ is modelled on Banach spaces of vector
bundle sections. The fibers of $E$ are manifolds diffeomorphic to $F$ that
can be endowed with Riemannian metrics coming from the metric in $F$. The
tangent spaces to the fibers form a vector bundle which we denote by $T^{f}E$%
. For $\xi \in E$ we write $\exp _{\xi }:T_{\xi }^{f}E\rightarrow E_{\xi }$
for the exponential of the Riemannian metric in the fiber $E_{\xi }$ through 
$\xi $.

Now, take a section $\sigma \in \Gamma E$ and write $\left( T^{f}E\right)
_{\sigma }$ for the restriction of $T^{f}E$ to $\mathrm{im}\sigma $. Then we
can define a map $\mathrm{Exp}_{\sigma }:\Gamma \left( T^{f}E\right)
_{\sigma }\rightarrow \Gamma E$ by associating to a section $Y:\mathrm{im}%
\sigma \rightarrow T^{f}E$ of $\left( T^{f}E\right) _{\sigma }$ the section 
\begin{equation*}
\mathrm{Exp}_{\sigma }\left( \tau \right) \left( x\right) =\exp _{\sigma
\left( x\right) }\left( Y\left( \sigma \left( x\right) \right) \right)
\end{equation*}%
of $E$. It turns out that $\mathrm{Exp}_{\sigma }$ is continuous and there
are neighborhoods $\mathcal{V}_{\sigma }\subset \Gamma \left( T^{f}E\right)
_{\sigma }$ of the zero section and $\mathcal{U}_{\sigma }\subset \Gamma E$
such that $\mathrm{Exp}_{\sigma }:\mathcal{V}_{\sigma }\rightarrow \mathcal{U%
}_{\sigma }$ is a homeomorphism.

Thus $\mathrm{Exp}_{\sigma }$ defines a chart for $\Gamma E$ around $\sigma $
modelled on the open set $\mathcal{V}$ of the Banach space $\Gamma \left(
T^{f}E\right) _{\sigma }$. It can be proved that the charts $\mathrm{Exp}%
_{\sigma }$, $\sigma \in \Gamma E$, form a differentiable atlas, so that $%
\Gamma E$ becomes a smooth manifold.

This smooth structure can be read off from the equivariant maps $f_{\sigma } 
$ of the sections $\sigma \in \Gamma E$. In fact, one can build on the space
of continuous maps $\mathcal{C}\left( R,F\right) $ the structure of a
differentiable manifold. It can be proved that the map 
\begin{equation*}
\sigma \in \Gamma E\mapsto \mathcal{C}_{\mathrm{eq}}\left( R,F\right)
\subset \mathcal{C}\left( R,F\right)
\end{equation*}%
is an embedding, so that $\mathcal{C}_{\mathrm{eq}}\left( R,F\right) $
becomes a submanifold diffeomorphic to $\Gamma E$.

A continuous automorphism $\phi \in \mathrm{Aut}\left( Q\right) $ lifts to a
map in $\Gamma \left( R\times _{K}F\right)$ defined by
\begin{equation*}
\Gamma \phi :\Gamma \left( R\times _{K}F\right) \rightarrow \Gamma
\left( R\times _{K}F\right), \hspace{10pt}
\Gamma \phi \left( \sigma \right) \left( x\right) =\phi \left( \sigma \left(
\phi ^{-1}\left( x\right) \right) \right) .
\end{equation*}%
The equivariant map of $\Gamma \phi \left( \sigma \right) $ is $f_{\sigma
}\circ \phi $ where $f_{\sigma }:Q\rightarrow F$ is the equivariant map of $%
\phi $. Since the map $f\in \mathcal{C}\left( Q,F\right) \mapsto f\circ \phi
\in \mathcal{C}\left( Q,F\right) $ is smooth it follows that these liftings
are smooth as well.

\begin{proposicao}
\label{propdifonsect}If $\phi \in \mathrm{Aut}\left( Q\right) $ then $\Gamma
\phi :\Gamma \left( R\times _{K}F\right) \rightarrow \Gamma \left( R\times
_{K}F\right) $ is differentiable with the above manifold structure.
\end{proposicao}

The differencial $d\left( \Gamma \phi \right) _{\sigma }$ at the section $%
\sigma $ is a linear map defined on the tangent space $T_{\sigma }\Gamma E$
which is the space of sections $\Gamma T_{\sigma }^{f}E$. It is given
pointwise as follows. The differentials $d\phi $ of $\phi $ along the fibers
are well defined giving rise to a bundle map $d\phi :T^{f}E\rightarrow
T^{f}E $. This bundle map extends to a linear map on the space of sections $%
\Gamma T^{f}E$, which we denote by $\Phi :\Gamma T^{f}E\rightarrow \Gamma
T^{f}E$. Then the following expression for $d\left( \Gamma \phi \right)
_{\sigma }$ can be proved.

\begin{proposicao}
\label{propdiffi}$d\left( \Gamma \phi \right) _{\sigma }$ is the restriction
of $\Phi $ to $T_{\sigma }^{f}E$.
\end{proposicao}

We turn now to the gauge group $\mathcal{G}=\mathcal{G}\left( Q\right) $ of $%
Q$ which is the (normal) subgroup of the group $\mathrm{Aut}\left( Q\right) $
whose elements are continuous automorphisms that induce the identity map on $%
X$. An element $\gamma \in \mathcal{G}$ defines the continuous map $%
f_{\gamma }:Q\rightarrow G$ by $\gamma \left( q\right) =q\cdot f_{\gamma
}\left( q\right) $, which is equivariant $f_{\gamma }\left( qg\right)
=g^{-1}f\left( q\right) g$. The equivariance property ensures that $%
f_{\gamma }$ is completely determined by its values in $R$. Now we endow $%
\mathcal{G}$ with the distance 
\begin{equation*}
d\left( \gamma ,\eta \right) =\sup_{q\in R}d_{0}\left( f_{\gamma }\left(
q\right) ,f_{\eta }\left( q\right) \right)
\end{equation*}%
where $d_{0}$ is a distance in $G$ coming from an invariant Riemannian
metric. $\mathcal{G}$ has a structure of Banach-Lie group compatible with
the topology given by $d$. Its Lie algebra $\mathrm{Lie}\left( \mathcal{G}%
\right) $ is the space $\Gamma \left( Q\times _{\mathrm{Ad}}\mathfrak{g}%
\right) $ of continuous sections of the vector bundle $Q\times _{\mathrm{Ad}}%
\mathfrak{g}$, which identifies with the space of continuous maps $%
f:Q\rightarrow \mathfrak{g}$ that are equivariant, $f\left( qg\right) =%
\mathrm{Ad}\left( g^{-1}\right) f\left( q\right) $. The exponential map in $%
\mathcal{G}$ is obtained pointwise by associating to $f:Q\rightarrow 
\mathfrak{g}$ the element $\gamma =\mathrm{Exp}f\in \mathcal{G}$ whose
equivariant map is $\exp f\left( \cdot \right) $, that is, $\gamma \left(
q\right) =q\exp f\left( q\right) $. The exponential map yields a coordinate
system $\mathrm{Exp}:\mathcal{V}\rightarrow \mathcal{U}$ for suitable chosen
open sets $\mathcal{V}\subset \mathrm{Lie}\left( \mathcal{G}\right) $ and $%
\mathcal{U}\subset \mathcal{G}$ containing the zero section and the identity
of $\mathcal{G}$, respectively.

The action of $\mathcal{G}\times \Gamma \left( Q\times _{G}F\right)
\rightarrow \Gamma \left( Q\times _{G}F\right) $ on a space of sections $%
\Gamma \left( Q\times _{G}F\right) $ is smooth. It induces an infinitesimal
action which associates to an element $Y\in \mathrm{Lie}\left( \mathcal{G}%
\right) $ a vector field $\widetilde{Y}$ on $\Gamma \left( Q\times
_{G}F\right) $ defined by 
\begin{equation*}
\widetilde{Y}\left( \sigma \right) =\frac{d}{dt}\left( \mathrm{Exp}tY\right)
\left( \sigma \right) _{\left\vert t=0\right. }.
\end{equation*}%
The map $Y\mapsto \widetilde{Y}\left( \sigma \right) $ is the differential
at $1$ of the evaluation map $E\left( \gamma \right) =\gamma \left( \sigma
\right) $.

If $\phi \in \mathrm{Aut}\left( Q\right) $ then the conjugation $C_{\phi }:%
\mathcal{G}\rightarrow \mathcal{G}$, $C_{\phi }\left( \gamma \right) =\phi
\gamma \phi ^{-1}$ is a smooth map of $\mathcal{G}$. In fact, if \ $%
f_{\gamma }$ is the equivariant map of $\gamma \in \mathcal{G}$ then $\phi
\gamma \phi ^{-1}\left( q\right) =q\cdot f_{\gamma }\left( \phi ^{-1}\left(
q\right) \right) $. Hence in a chart $\mathrm{Exp}:\mathcal{V}\rightarrow 
\mathcal{U}$ around the identity $C_{\phi }$ becomes the map that associates 
$f:Q\rightarrow \mathfrak{g}$ to $f\circ \phi ^{-1}$. If we view $\mathrm{Lie%
}\left( \mathcal{G}\right) $ as the space sections $\Gamma \left( Q\times _{%
\mathrm{Ad}}\mathfrak{g}\right) $ then the conjugation $C_{\phi }$ is read
in a chart as the restriction of a bounded linear map, namely $\phi \circ
\sigma \circ \phi ^{-1}$, $\sigma \in \Gamma \left( Q\times _{\mathrm{Ad}}%
\mathfrak{g}\right) $. Therefore, $C_{\phi }$ is smooth at the identity of $%
\mathcal{G}$. Since it is an automorphism of $\mathcal{G}$ it is smooth
everywhere. 

\section{Perturbation of the attractor Morse component}

As mentioned in the introduction by the results of \cite{smbflow} and \cite%
{msm} show that in a flag bundle $\mathbb{E}_{\Theta }$ there is a finest
Morse decomposition, which is given by a continuous section $\chi _{\phi }$
of the associated bundle $Q\times _{G}\mathbb{A}\left( H_{\mathrm{Mo}}\left(
\phi \right) \right) $ where $H_{\mathrm{Mo}}\left( \phi \right) \in \mathrm{%
cl}\mathfrak{a}^{+}$ and $\mathbb{A}\left( H_{\mathrm{Mo}}\left( \phi
\right) \right) $ is its adjoint orbit. The flag type of $\phi $ is the flag
manifold $\mathbb{F}_{\Theta _{\mathrm{Mo}}\left( \phi \right) }$ where the
subset $\Theta _{\mathrm{Mo}}\left( \phi \right) $ of the simple system of
roots $\Sigma $ is given by 
\begin{equation*}
\Theta _{\mathrm{Mo}}\left( \phi \right) =\{\alpha \in \Sigma :\alpha \left(
H_{\mathrm{Mo}}\left( \phi \right) \right) =0\}.
\end{equation*}%
Since the components of the finest Morse decomposition are constructed by
the fixed point sets of conjugates of $h=\exp H_{\mathrm{Mo}}\left( \phi
\right) $, it follows that in the specific flag bundle $\mathbb{E}_{\Theta _{%
\mathrm{Mo}}\left( \phi \right) }\rightarrow X$ the attractor Morse
component meets each fiber in a singleton.

Hence this attractor component is the image of a continuous section $\sigma
_{\mathrm{Mo}}\in \Gamma \mathbb{E}_{\Theta _{\mathrm{Mo}}\left( \phi
\right) }$, which is $\phi $-invariant. This means that $\sigma _{\mathrm{Mo}%
}$ is a fixed point for the action of $\phi $ on $\Gamma \mathbb{E}_{\Theta
_{\mathrm{Mo}}\left( \phi \right) }$.

The existence of fixed points on $\Gamma \mathbb{E}_{\Theta _{\mathrm{Mo}%
}\left( \phi \right) }$ for automorphisms close to $\phi $ will be ensured
by the following easy consequence of the Implicit Function Theorem.

\begin{lema}
\label{lemteoimplic}Let $M$ and $\mathcal{H}$ be $\mathcal{C}^{k}$, $k\geq 1$%
, Banach manifolds and $p:\mathcal{H}\times M\rightarrow M$ a $C^{k}$ map.
For $y\in \mathcal{H}$ let $\gamma _{y}:M\rightarrow M$ be the partial map $%
\gamma _{y}\left( x\right) =p\left( y,x\right) $. Take $\left(
y_{0},x_{0}\right) \in \mathcal{H}\times M$ such that $\gamma _{y_{0}}\left(
x_{0}\right) =x_{0}$ and assume that $d\gamma _{y_{0}}-\mathrm{id}$ is
invertible. Then there exists a neighborhood $V$ of $y_{0}$ in $\mathcal{H}$
and a $\mathcal{C}^{k}$ map $\chi :V\rightarrow M$ such that $\chi \left(
y\right) $ is fixed by $\gamma _{y}$, $y\in V$, that is, $\gamma _{y}\left(
\chi \left( y\right) \right) =\chi \left( y\right) $.
\end{lema}

\begin{profe}
Take local coordinates around $y_{0}$ and $x_{0}$ and apply the implicit
function theorem to the map%
\begin{equation*}
f\left( \gamma ,x\right) =\gamma \left( x\right) -x,
\end{equation*}
which is possible by the existence of $\left( d\gamma _{y_{0}}-\mathrm{id}%
\right) ^{-1}$.
\end{profe}

We will apply Lemma \ref{lemteoimplic} to the action of $\mathcal{G}\phi
\times \Gamma \mathbb{E}_{\Theta _{\mathrm{Mo}}}\rightarrow \Gamma \mathbb{E}%
_{\Theta _{\mathrm{Mo}}}$, which by Proposition \ref{propdifonsect} is
smooth. We take $y_{0}=\phi $ and $x_{0}$ to be the section $\sigma _{%
\mathrm{Mo}}\in \Gamma \mathbb{E}_{\Theta _{\mathrm{Mo}}}$ such that $%
\mathcal{M}_{\Theta _{\mathrm{Mo}}}^{+}=\mathrm{im}\sigma _{\mathrm{Mo}}$ is
the attractor Morse component in $\mathbb{E}_{\Theta _{\mathrm{Mo}}}$. To
check the invertibility condition of Lemma \ref{lemteoimplic} recall that by
Proposition \ref{propdiffi} the differential $d\left( \Gamma \phi \right)
_{\sigma }$ is the restriction of $\Phi $ to $T_{\sigma }^{f}E$. The
invertibility of $\Phi -\mathrm{id}$ is a consequence of the following
uniformity lemma.

\begin{lema}
\label{lemexistnorma}There exists a norm $\left\Vert \cdot \right\Vert $ in $%
T^{f}\mathbb{E}_{\Theta }$ and $\alpha >0$ such that 
\begin{equation*}
\left\Vert \Phi \right\Vert =\sup_{\eta \in \mathcal{M}_{\Theta _{\mathrm{Mo}%
}}^{+}}\left\Vert d\phi _{\eta }\right\Vert \leq e^{-\alpha }<1.
\end{equation*}
\end{lema}

\begin{profe}
It was proved in \cite{conl}, Corolary 5.2 that there are a norm $\left\Vert
\cdot \right\Vert $ in $T^{f}\mathbb{E}_{\Theta _{\mathrm{Mo}}}$ and $\alpha
>0$ such that if $v\in T_{\eta }^{f}\mathbb{E}_{\Theta _{\mathrm{Mo}}}$ with 
$\eta \in \mathcal{M}_{\Theta _{\mathrm{Mo}}}^{+}$ then 
\begin{equation*}
\left\Vert d\phi _{\eta }^{n}v\right\Vert <e^{-\alpha n}\left\Vert
v\right\Vert .
\end{equation*}%
In particular $\left\Vert d\phi _{\eta }v\right\Vert <e^{-\alpha }\left\Vert
v\right\Vert $ so that for any $\eta \in \mathcal{M}_{\Theta _{\mathrm{Mo}%
}}^{+}$ we have $\left\Vert d\phi _{\eta }\right\Vert <e^{-\alpha }<1$. Thus 
$\left\Vert \Phi \right\Vert \leq e^{-\alpha }<1$ as claimed.
\end{profe}

\begin{corolario}
$d\left( \Gamma \phi \right) _{\sigma }-\mathrm{id}$ is invertible.
\end{corolario}

Now we can apply the Implicit Function Theorem via Lemma \ref{lemteoimplic}
to get fixed points of $\gamma \phi $.

\begin{proposicao}
\label{propptofix}With the above notation, there exists a neighborhood $V$
of the identity of $\mathcal{G}$ and a smooth map $s:V\rightarrow \Gamma 
\mathbb{E}_{\Theta _{\mathrm{Mo}}}$, $s\left( \gamma \right) =\sigma
_{\gamma }$, such that $\sigma _{\gamma }$ is $\gamma \phi $-invariant and a 
$\Gamma \left( \gamma \phi \right) $-fixed point. \ Its differential at the
identity $1\in \mathcal{G}$ is%
\begin{equation*}
ds_{1}\left( Y\right) =-\left( \Phi -\mathrm{id}\right) ^{-1}\widetilde{Y}%
\left( \sigma _{1}\right)
\end{equation*}%
where $Y\in \mathrm{Lie}\left( \mathcal{G}\right) $ and $\widetilde{Y}\left(
\tau \right) =d/dt\left( \mathrm{Exp}Y\right) \left( \tau \right) _{t=0}$.
\end{proposicao}

\begin{profe}
Apply Lemma \ref{lemteoimplic} to the map $\mathcal{G}\phi \times \Gamma
\left( \mathbb{E}_{\Theta _{\mathrm{Mo}}}\right) \rightarrow \Gamma \left( 
\mathbb{E}_{\Theta _{\mathrm{Mo}}}\right) $ given by the action of $\mathcal{%
G}\phi $ on $\Gamma \left( \mathbb{E}_{\Theta _{\mathrm{Mo}}}\right) $,
which is possible because $d\left( \Gamma \phi \right) _{\sigma }-\mathrm{id}
$ is invertible.
\end{profe}

We proceed now to prove that an invariant section $\sigma _{\gamma }$ of
Proposition \ref{propptofix} is the attractor Morse component of $\gamma
\phi $. This fact will be used below to give a differentiable
parametrization of the invariant measures for the flows generated by $\gamma
\phi $.

The proof that $\sigma _{\gamma }$ is the attractor Morse component of $%
\gamma \phi $ uses the concept of transversality of sections in flag
bundles, which can be seen, for instance in \cite{alvsm}, Section 6.2.

Take $\psi \in \mathrm{Aut}\left( Q\right) $ a flag bundle $\mathbb{E}%
_{\Theta }$ and its dual $\mathbb{E}_{\Theta ^{\ast }}$. Suppose that threre
are transversal $\psi $-invariant sections $\sigma :X\rightarrow \mathbb{E}%
_{\Theta }$ and $\sigma ^{\ast }:X\rightarrow \mathbb{E}_{\Theta ^{\ast }}$.
Take $H_{\Theta }\in \mathrm{cl}\mathfrak{a}^{+}$ such that 
\begin{equation*}
\Theta =\{\alpha \in \Sigma :\alpha \left( H_{\Theta }\right) =0\}
\end{equation*}%
and let $Z_{H_{\Theta }}$ stand for the centralizer of $H_{\Theta }$ in $G$.

The transversal sections $\sigma $ and $\sigma ^{\ast }$ give rise to a $%
\psi $-invariant section $\chi $ of the bundle $Q\times _{G}$\textrm{Ad}$%
\left( G\right) H_{\Theta }$, because \textrm{Ad}$\left( G\right) H_{\Theta
} $ \ identifies to the set of transversal elements in $\mathbb{F}_{\Theta
}\times \mathbb{F}_{\Theta ^{\ast }}$. In turn $\chi $ defines a $%
Z_{H_{\Theta }}$-reduction $L\subset Q$, thas is $\psi $-invariant (see \cite%
{conl}, \cite{smsec}).

The $Z_{H}$-reduction yields a conjugation around $\sigma $ between $\psi $
and its linearization $\Psi =d\psi $ on the vertical tangent bundle $%
T_{\sigma }^{f}\mathbb{E}_{\Theta }$ (see \cite{conl}):

\begin{lema}
The linearization $\Psi $ is conjugate to $\psi $, that is, there exists a
homeomorphism $e:T_{\sigma }^{f}\mathbb{E}_{\Theta }\rightarrow S$ where $%
S\subset \mathbb{E}_{\Theta }$ is open such that $e\Psi \left( v\right)
=\psi \left( ev\right) $ if $v\in T_{\sigma }^{f}\mathbb{E}_{\Theta }$. The
conjugation $e$ maps the zero section $0$ to itself. The same holds for the
section $\sigma ^{\ast }$ of $\mathbb{E}_{\Theta ^{\ast }}$.
\end{lema}

By using this conjugation we can prove that the image of $\sigma $ is the
attractor Morse component of $\psi $ in $\mathbb{E}_{\Theta }$.

\begin{proposicao}
\label{propsecisMo}Let $\psi \in \mathrm{Aut}\left( Q\right) $, $\mathbb{E}%
_{\Theta }$ a flag bundle with dual $\mathbb{E}_{\Theta ^{\ast }}$. Assume
the following conditions:

\begin{enumerate}
\item There are transversal $\psi $-invariant sections $\sigma :X\rightarrow 
\mathbb{E}_{\Theta }$ and $\sigma ^{\ast }:X\rightarrow \mathbb{E}_{\Theta
^{\ast }}$.

\item $\left\Vert d\psi _{\sigma }\right\Vert <1$, where $d\psi _{\sigma
}=\Psi $ is seen as a linear map in the space of sections of the vertical
tangent bundle at the image of $\sigma $.
\end{enumerate}

Then, $\mathrm{im}\sigma $ is the attractor Morse component of the finest
Morse decomposition of $\psi $ in $\mathbb{E}_{\Theta }$.
\end{proposicao}

\begin{profe}
Let $B\subset T_{\sigma }^{f}\mathbb{E}_{\Theta }$ be the unit ball and $%
A=e\left( B\right) $, where $e$ is the conjugation of the above lemma.
Assumption (2) implies that $A$ is an attractor neighborhood of $\mathrm{im}%
\sigma $. In fact, if $x_{n}=e\left( v_{n}\right) \in A$ then 
\begin{equation*}
\left\Vert \Psi ^{n}v_{n}\right\Vert \leq \left\Vert \Psi ^{n}\right\Vert
\cdot \left\Vert v_{n}\right\Vert \leq \left\Vert \Psi \right\Vert
^{n}\rightarrow 0.
\end{equation*}%
Hence the accumulation points of $\psi \left( x_{n}\right) =e\Psi \left(
v_{n}\right) $ are contained in $\mathrm{im}\sigma $. It follows that $%
\mathrm{im}\sigma $ contains the attractor component $\mathcal{M}_{\Theta
}^{+}\left( \psi \right) $ of the finest Morse in $\mathbb{E}_{\Theta }$,
which implies that $\mathrm{im}\sigma =\mathcal{M}_{\Theta }^{+}\left( \psi
\right) $ because $\mathrm{im}\sigma $ meets each fiber in a sigleton.
\end{profe}

\noindent
\textbf{Remark:} It follows by the above proposition that the flag type $%
\Theta _{\mathrm{Mo}}\left( \psi \right) $ of $\psi $ is contained in $%
\Theta $, since its attractor component $\mathcal{M}_{\Theta }^{+}\left(
\psi \right) $ is the image of a section.

Now we can check that the sections of Proposition \ref{propptofix} are
indeed Morse components.

\begin{teorema}
\label{teosecMorse}Let $\Theta _{\mathrm{Mo}}=\Theta _{\mathrm{Mo}}\left(
\phi \right) $ be the flag type of $\phi $ and $\sigma \in \Gamma \left( 
\mathbb{E}_{\Theta _{\mathrm{Mo}}}\right) $ the section defining the
attractor Morse component of $\phi $ in $\mathbb{E}_{\Theta _{\mathrm{Mo}}}$%
. Then we can choose $V\subset \mathcal{G}$ for the map $\gamma \in V\mapsto
\sigma _{\gamma }\in \Gamma \left( \mathbb{E}_{\Theta _{\mathrm{Mo}}}\right) 
$ of Proposition \ref{propptofix}, such that the attractor Morse component
of $\gamma \phi $ in $\mathbb{E}_{\Theta _{\mathrm{Mo}}}$ is the image of $%
\sigma _{\gamma }$, $\gamma \in V$.
\end{teorema}

\begin{profe}
Is a consequence of Proposition \ref{propsecisMo}. To check its assumptions
we note first that the flag type of $\phi ^{-1}$ is $\Theta _{\mathrm{Mo}%
}^{\ast }$. Hence there exists a section $\sigma _{\mathrm{Mo}}^{\ast
}:X\rightarrow \mathbb{E}_{\Theta _{\mathrm{Mo}}^{\ast }}$ which is $\phi
^{-1}$-invariant and its image is the attractor Morse component of $\phi
^{-1}$ in $\mathbb{E}_{\Theta _{\mathrm{Mo}}^{\ast }}$. The section $\sigma
_{\mathrm{Mo}}^{\ast }$ is transversal to $\sigma _{\mathrm{Mo}}$.

Now we can apply Proposition \ref{propptofix} to $\phi $ and $\phi ^{-1}$
and the sections $\sigma _{\mathrm{Mo}}$ and $\sigma _{\mathrm{Mo}}^{\ast }$%
, respectively, to get a neighborhood of the identity $W\subset \mathcal{G}$
and differentiable maps $\gamma \in W\mapsto \sigma _{\gamma }\in \Gamma
\left( \mathbb{E}_{\Theta _{\mathrm{Mo}}}\right) $ and $\gamma \in W\mapsto
\sigma _{\gamma }^{\ast }\in \Gamma \left( \mathbb{E}_{\Theta _{\mathrm{Mo}%
}^{\ast }}\right) $ such that $\sigma _{\gamma }$ is $\gamma \phi $%
-invariant and $\sigma _{\gamma }^{\ast }$ is $\gamma ^{-1}\phi ^{-1}$%
-invariant. Since transversality is an open condition and $\sigma _{\mathrm{%
Mo}}^{\ast }$ is transversal to $\sigma _{\mathrm{Mo}}$ we can shrink $W$ so
that $\sigma _{\gamma }$ and $\sigma _{\eta }^{\ast }$ are transversal if $%
\gamma ,\eta \in W$.

Put $V=W\cap \phi W\phi ^{-1}$. If $\gamma \in V$ then $\sigma _{\gamma }$
is $\gamma \phi $-invariant. The section $\sigma _{\eta }^{\ast }\in \Gamma
\left( \mathbb{E}_{\Theta _{\mathrm{Mo}}^{\ast }}\right) $ is also $\gamma
\phi $-invariant where $\eta =\phi \gamma \phi ^{-1}$, because $\sigma
_{\eta }^{\ast }$ is invariant by $\eta ^{-1}\phi ^{-1}$ and hence by $\phi
\eta =\left( \eta ^{-1}\phi ^{-1}\right) ^{-1}$. So $\sigma _{\eta }^{\ast }$
is invariant by $\gamma \phi =\left( \phi \eta \phi ^{-1}\right) \phi =\phi
\eta $. Hence the first assumption of Proposition \ref{propsecisMo} holds
for $\psi =\gamma \phi $ if $\gamma \in V$.

Now, by continuity we can take $V$ so that $\left\Vert d\left( \gamma \phi
\right) _{\sigma _{\gamma }}\right\Vert <1$ (see Lemma \ref{lemexistnorma}).
Hence the second assumption of Proposition \ref{propsecisMo} holds so that $%
\mathrm{im}\sigma _{\gamma }$ is the attractor Morse component in $\mathbb{E}%
_{\Theta }$ of $\gamma \phi $, $\gamma \in V$.
\end{profe}

\vspace{12pt}%

\noindent%
\textbf{Remark:} It follows by the previous theorem that the flag type $%
\Theta _{\mathrm{Mo}}\left( \gamma \phi \right) $ is contained in $\Theta _{%
\mathrm{Mo}}\left( \phi \right) $ if $\gamma $ is close to identity. This
fact however can be obtained by general topological facts about Morse
decompositions.

\section{Smoothness of the additive cocycles\label{secociclan}}

The objective here is to prove that the additive cocycles $\mathsf{a}%
_{\omega }^{\phi }\left( n,\cdot \right) $ on a flag bundle $\mathbb{E}%
_{\Theta }$ depend differentiably on the flow $\phi $ in the following
sense: For each $\gamma \in \mathcal{G}$ we have the cocycle $\mathsf{a}%
_{\omega }^{\gamma \phi }\left( n,\xi \right) $ on $\mathbb{E}_{\Theta }$
over the flow induced by $\gamma \phi $, where $\omega $ belongs to the
linear span of $\Omega \setminus \Omega _{\Theta }$. The map $\xi \in 
\mathbb{E}_{\Theta }\mapsto \mathsf{a}_{\omega }^{\gamma \phi }\left( n,\xi
\right) \in \mathbb{R}$ is continuous. Taking $n=1$ and fixing $\phi $ we
have a map 
\begin{equation}
F_{\omega }^{\phi }:\mathcal{G}\rightarrow \mathcal{C}\left( \mathbb{E}%
_{\Theta },\mathbb{R}\right) \qquad F_{\omega }^{\phi }\left( \gamma \right)
=\mathsf{a}_{\omega }^{\gamma \phi }\left( 1,\cdot \right) ,
\label{fordefFtetafi}
\end{equation}%
where $\mathcal{C}\left( \mathbb{E}_{\Theta },\mathbb{R}\right) $ is the
Banach space of continuous functions $\mathbb{E}_{\Theta }\rightarrow 
\mathbb{R}$ endowed with the $\sup $ norm.

\begin{proposicao}
\label{propdifcocycle}The map $F_{\omega }^{\phi }:\mathcal{G}\rightarrow 
\mathcal{C}\left( \mathbb{E}_{\Theta },\mathbb{R}\right) $ is smooth. Its
differential at $\gamma =1$ is 
\begin{equation}
d\left( F_{\omega }^{\phi }\right) _{1}\left( Y\right) \left( r\cdot
b_{\Theta }\right) =\omega \circ \mathbf{a}\left( f_{Y}\left( \phi
^{R}\left( r\right) \right) \right)  \label{fordiferena}
\end{equation}%
where $Y\in T_{1}\mathcal{G}$ is a section of the bundle $Q\times _{\mathrm{%
Ad}}\mathfrak{g}$, $f_{Y}:R\rightarrow \mathfrak{g}$ is the equivariant map
of $Y$ and $\mathbf{a}:\mathfrak{g}\rightarrow \mathfrak{a}$ is the
projection onto the $\mathfrak{a}$-component of the Iwasawa decomposition $%
\mathfrak{g}=\mathfrak{k}\oplus \mathfrak{a}\oplus \mathfrak{n}$. Here $%
b_{\Theta }$ is the origin of $\mathbb{F}_{\Theta }$.
\end{proposicao}

To prove the proposition consider first the map $F^{\phi }$ defined by the
cocycle over $\phi ^{R}$ on the reduction $R$. Let $i:\mathcal{G}\rightarrow 
\mathcal{C}\left( R,G\right) $ be the embedding $i\left( \gamma \right)
=f_{\gamma }$ where $f_{\gamma }$ is the equivariant map defined by $\gamma
(r)=rf_{\gamma }(r)$, $r\in R$.

Given the Iwasawa decomposition $G=KAN$ let $\mathsf{A}:G\rightarrow A$ be
the projection, which is a smooth map, and define the map $\mathcal{A}:%
\mathcal{C}(R,G)\rightarrow \mathcal{C}(R,A)$, $\mathcal{A}(f)=\mathsf{A}%
\circ f$, which is smooth as well.

By definition of the cocycle $\mathsf{a}^{\gamma }(n,r)$ over $R\subset Q$,
defined by $\gamma \in \mathcal{G}$, we have 
\begin{equation*}
\mathsf{a}^{\gamma }(1,r)=\log \left( \mathcal{A}(i(\gamma ))(r)\right)
\qquad r\in R.
\end{equation*}%
Hence the map 
\begin{equation}
\begin{array}{cccl}
F: & \mathcal{G} & \longrightarrow & \mathcal{C}(R,\mathfrak{a}) \\ 
& \gamma & \longmapsto & \mathsf{a}^{\gamma }(1,\cdot )%
\end{array}
\label{forZ}
\end{equation}%
is smooth.

As to the map $F^{\phi }(\gamma )=\mathsf{a}^{\gamma \phi }(1,\cdot )$, we
have 
\begin{equation}
\mathsf{a}^{\gamma \phi }(1,r)=\mathsf{a}^{\gamma }\left( 1,\phi
^{R}(r)\right) +\mathsf{a}^{\phi }(1,r),\qquad r\in R.  \label{fora}
\end{equation}%
In fact $\phi (r)=\phi ^{R}(r)\cdot hn$, with $hn\in AN$ and $\log h=\mathsf{%
a}^{\phi }(1,r)$. Also, 
\begin{equation*}
\gamma \left( \phi ^{R}(r)\right) =\phi ^{R}(r)\cdot f_{\gamma }\left( \phi
^{R}(r)\right) .
\end{equation*}%
Write the Iwasawa decomposition of $f_{\gamma }\left( \phi ^{R}(r)\right)
\in G$ as $f_{\gamma }\left( \phi ^{R}(r)\right) =kan^{\prime }$, so that $%
\log a=\mathsf{a}^{\gamma }\left( 1,\phi ^{R}(r)\right) $. Then we get 
\begin{equation*}
\gamma \phi (r)=\gamma \left( \phi ^{R}(r)\cdot hn\right) =\gamma \left(
\phi ^{R}(r)\right) \cdot hn=\phi ^{R}(r)\cdot kah\overline{n}
\end{equation*}%
with $\overline{n}\in N$. Hence, $\mathsf{a}^{\gamma \phi }(1,r)=\log a+\log
h=\mathsf{a}^{\gamma }\left( 1,\phi ^{R}(r)\right) +\mathsf{a}^{\phi }(1,r)$
as claimed.

By (\ref{fora}) the map $F^{\phi }$ is the composition of $F$ with the
affine map of $\mathcal{C}(R,\mathfrak{a})$ given by $f\mapsto f\circ \phi
^{R}+\mathsf{a}^{\phi }(1,\cdot )$. Its linear part $f\in \mathcal{C}(R,%
\mathfrak{a})\mapsto f\circ \phi ^{R}\in \mathcal{C}(R,\mathfrak{a})$ is an
isometry hence the affine map is smooth showing that $F^{\phi }$ is smooth.

The composition $\omega \circ F^{\phi }$, $\omega \in \mathrm{span}\left(
\Omega \setminus \Omega _{\Theta }\right) $ is also smooth and since factors
to $\mathbb{E}_{\Theta }$ yielding $F_{\omega }^{\phi }$ we conclude that
this last map is smooth as well.

Now we look at the differential of $F_{\omega }^{\phi }$. We have 
\begin{equation*}
d\left( F_{\omega }^{\phi }\right) _{1}\left( Y\right) =\frac{d}{dt}%
F_{\omega }^{\phi }\left( \exp tY\right) _{t=0}
\end{equation*}
where $\exp tY$ is the exponential in $\mathcal{G}$, whose equivariant
function is $f_{\exp tY}\left( r\right)=e^{tf_{Y}\left(r\right)}$. Since
we know in advance that $F_{\omega }^{\phi }$ is differentiable we can
perform the derivative in the right hand side pointwise. By (\ref{fora}) 
\begin{equation*}
F_{\omega }^{\phi }\left( \exp tY\right) \left( r\cdot b_{\Theta }\right)
=\omega \circ \mathsf{a}^{\gamma _{t}}\left( 1,\phi ^{R}(r)\cdot b_{\Theta
}\right) +\omega \circ \mathsf{a}^{\phi }(1,r\cdot b_{\Theta })
\end{equation*}%
where $\gamma _{t}=\exp tY$. But $\mathsf{a}^{\gamma _{t}}\left( 1,\phi
^{R}(r)\right) =\log h_{t}$ where $\gamma _{t}\left( \phi ^{R}(r)\right)
=\phi ^{R}\left( r\right) \cdot u_{t}h_{t}n_{t}$ is the Iwasawa
decomposition. But 
\begin{eqnarray*}
\gamma _{t}\left( \phi ^{R}(r)\right) &=&\phi ^{R}(r)\cdot f_{\exp tY}\left(
\phi ^{R}\left( r\right) \right) \\
&=&\phi ^{R}(r)\cdot e^{tf_{Y}\left( \phi ^{R}\left( r\right) \right) }
\end{eqnarray*}%
Hence $u_{t}h_{t}n_{t}=e^{tf_{Y}\left( \phi ^{R}\left( r\right) \right) }$.
Taking derivatives it follows that $\frac{d}{dt}\mathsf{a}^{\gamma
_{t}}\left( 1,\phi ^{R}(r)\right) _{t=0}$ is the $\mathfrak{a}$-component of 
$f_{Y}\left( \phi ^{R}\left( r\right) \right) $.

\section{Main theorem}

In this section we prove the main result on the differentiability of
Lyapunov spectra. As before we let $\nu $ be a $\phi $-invariant probability
measure on $X$ and write $H_{\phi }^{+}\left( x\right) \in \mathrm{cl}%
\mathfrak{a}^{+}$ for the polar exponent of $\phi \in \mathrm{Aut}\left(
Q\right) $ given by the subadditive ergodic theorem. Also $\Omega _{\Theta }$
is the set of fundamental weights corresponding to the simple roots in $%
\Theta $, so that $\mathrm{span}\left( \Omega \setminus \Omega _{\Theta
}\right) $ is the annihilator of $\mathfrak{a}\left( \Theta \right) $ (see (%
\ref{foradeteta})).

\begin{teorema}
\label{teoprinc}Let $\Theta =\Theta _{\mathrm{Mo}}\left( \phi \right) $ be
the flag type of the Morse decompositions of $\phi $. Take $\omega \in
\Omega \setminus \Omega _{\Theta }\mathfrak{\ }$. Then the map $\Lambda
_{\omega }:\mathcal{G}\rightarrow \mathbb{R}$, 
\begin{equation}
\Lambda _{\omega }\left( \gamma \right) =\int_{X}\omega \left( H_{\gamma
\phi }^{+}\left( x\right) \right) \nu \left( dx\right)  \label{forintegral}
\end{equation}%
is smooth at the identity $1\in \mathcal{G}$.
\end{teorema}

\vspace{12pt}%

\noindent%
\textbf{Remark:} The gauge group $\mathcal{G}$ of $X\times G$ is given
simply by maps $f:X\rightarrow G$ so that a $\gamma \in \mathcal{G}$ is
written $\gamma \left( x,g\right) =\left( x,gf_{\gamma }\left( x,g\right)
\right) $ where $f_{\gamma }\left( x,g\right) =g^{-1}f_{\gamma }\left(
x\right) g$ is the equivariant map and $f_{\gamma }\left( x\right)
=f_{\gamma }\left( x,1\right) $ is its restriction. The differentiable
structure on $\mathcal{G}$ is then constructed as a manifold of maps $%
\mathcal{C}\left( X,G\right) $. If $\phi $ is defined by the cocycle $\rho $
then $\gamma \phi \left( x,g\right) =\left( \tau \left( x\right) ,f_{\gamma
}\left( x\right) \rho \left( 1,x\right) g\right) $ so that $\gamma \phi $ is
defined by the cocycle generated by $f_{\gamma }\left( x\right) \rho \left(
1,x\right) $. Hence the map (\ref{forintegral}) is seen as a map $\mathcal{C}%
\left( X,G\right) \rightarrow \mathbb{R}$ and differentiability is taken
with respect to the manifold structure of $\mathcal{C}\left( X,G\right) $.

The first step in the proof of the theorem is to write the integral in (\ref%
{forintegral}) as an integral with respect to an invariant measure on the
flag bundle $\mathbb{E}_{\Theta _{\mathrm{Mo}}\left( \phi \right) }$. Recall
that by the multiplicative ergodic theorem version of \cite{alvsm} if $\psi
\in \mathrm{Aut}\left( Q\right) $ has invariant measure $\nu $ on $X $ then

\begin{enumerate}
\item for $\nu $-almost all $x$ there exists a partially defined Oseledets
section $\chi \left( x\right) \in Q\times _{G}\mathrm{Ad}\left( G\right)
H_{\psi }^{+}\left( x\right) $ such that the $\mathfrak{a}$-Lyapunov
exponent $\lambda _{\psi }\left( \xi \right) $ in the direction of $\xi $ in
the full flag bundle $\mathbb{E}$ is given by por 
\begin{equation}
\lambda _{\psi }\left( \xi \right) =\lim \frac{1}{n}\mathsf{a}^{\psi }\left(
n,\xi \right) =w^{-1}H_{\psi }^{+}\left( x\right) .  \label{forlyappolar}
\end{equation}%
if $\xi \in \mathrm{fix}\left( \chi \left( x\right) ,w\right) \subset 
\mathbb{E}$. Here $\mathsf{a}^{\psi }\left( n,\xi \right) $ is the cocycle
on $\mathbb{E}$ over $\psi $. In particular $H_{\psi }^{+}\left( x\right)
=\lambda _{\psi }\left( \xi \right) $ if $w=1$.

\item In this decomposition a component $\mathrm{fix}\left( \chi \left(
x\right) ,w\right) $ is contained in the Morse component $\mathcal{M}%
_{w}\left( \psi \right) $, so that the attractor Morse component $\mathcal{M}%
^{+}\left( \psi \right) $ contains $\mathrm{fix}\left( \chi \left( x\right)
,1\right) $.
\end{enumerate}

A similar picture holds in a partial flag bundle $\mathbb{E}_{\Theta }$. If $%
\omega \in \mathrm{span}\left( \Omega \setminus \Omega _{\Theta }\right) $
then the $\mathbb{R}$-valued cocycle $\mathsf{a}_{\omega }^{\psi }\left(
n,\xi \right) =\omega \left( \mathsf{a}^{\psi }\left( n,\xi \right) \right) $%
, \ factors to $\mathbb{E}_{\Theta }$ and we get 
\begin{equation*}
\lim \frac{1}{n}\mathsf{a}_{\omega }^{\psi }\left( n,\xi \right) =\omega
\left( w^{-1}H_{\psi }^{+}\left( x\right) \right)
\end{equation*}%
if $\xi \in \mathrm{fix}_{\Theta }\left( \chi \left( x\right) ,w\right)
\subset \mathbb{E}_{\Theta }$, which is the projection of $\mathrm{fix}%
\left( \chi \left( x\right) ,w\right) $.

Now, let $\Theta _{\mathrm{Mo}}\left( \psi \right) $ be the flag type of $%
\psi $ and take $\Theta \supset \Theta _{\mathrm{Mo}}\left( \psi \right) $.
In the flag bundle $\mathbb{E}_{\Theta }$ the attractor Morse component $%
\mathcal{M}_{\Theta }^{+}\left( \psi \right) $ is the image of a continuous
section $\sigma _{\Theta }^{\psi }:X\rightarrow \mathbb{E}_{\Theta }$ and
since $\mathrm{fix}_{\Theta }\left( \chi \left( x\right) ,1\right) \subset 
\mathcal{M}_{\Theta }^{+}\left( \psi \right) $ it follows that 
\begin{equation*}
\mathrm{fix}_{\Theta }\left( \chi \left( x\right) ,1\right) =\sigma _{\Theta
}^{\psi }\left( x\right)
\end{equation*}%
for any $x$ in the set of full measure where $\chi \left( x\right) $ is
defined. Hence if $\omega \in \mathrm{span}\left( \Omega \setminus \Omega
_{\Theta }\right) $ then for $\nu $-almost all $x$ we have 
\begin{equation*}
\lim \frac{1}{n}\mathsf{a}_{\omega }^{\psi }\left( n,\sigma _{\Theta }\left(
x\right) \right) =\omega \left( H_{\psi }^{+}\left( x\right) \right) .
\end{equation*}

Hence if we let $\nu _{\Theta }^{\psi }=\sigma _{\Theta }^{\psi }\cdot \nu $
be the push forward of $\nu $ under $\sigma _{\Theta }^{\psi }$ then the
Birkhoff ergodic theorem applied to $f\left( \xi \right) =\mathsf{a}_{\omega
}^{\psi }\left( 1,\xi \right) $ yields that 
\begin{equation}
\int \omega \left( H_{\psi }^{+}\left( x\right) \right) \nu \left( dx\right)
=\int \mathsf{a}_{\omega }^{\psi }\left( 1,\xi \right) \nu _{\Theta }^{\psi
}\left( d\xi \right) .  \label{forintegralflag}
\end{equation}

Summarizing we have

\begin{proposicao}
\label{propintflag}Let $\psi \in \mathrm{Aut}\left( Q\right) $ with
invariant measure $\nu $ on $X$. Take $\Theta \subset \Sigma $ with $\Theta
_{\mathrm{Mo}}\left( \psi \right) \subset \Theta $ so that the attractor
Morse component $\mathcal{M}_{\Theta }^{+}\left( \psi \right) $ in $\mathbb{E%
}_{\Theta }$ is the image of the continuous section $\sigma _{\Theta }^{\psi
}:X\rightarrow \mathbb{E}_{\Theta }$. Put $\nu _{\Theta }^{\psi }=\sigma
_{\Theta }^{\psi }\cdot \nu $ then (\ref{forintegralflag}) holds.
\end{proposicao}

\subsection{Smooth parametrization of the invariant measures in $\mathbb{E}%
_{\Theta }$}

Consider now the flag bundle $\mathbb{E}_{\Theta _{\mathrm{Mo}}}$ with $%
\Theta _{\mathrm{Mo}}=\Theta _{\mathrm{Mo}}\left( \phi \right) $ defining
the flag type of $\phi $. By Proposition \ref{propptofix} there exists a
neighborhood $V$ of the identity of $\mathcal{G}$ and a map $\gamma \in 
\mathcal{G}\mapsto \sigma _{\gamma }\in \Gamma \mathbb{E}_{\Theta _{\mathrm{%
Mo}}}$ such that the image of $\sigma _{\gamma }$ is the attractor Morse
component of $\gamma \phi $. Given the invariant measure $\nu $ on $X$ write 
$\nu _{\gamma }=\sigma _{\gamma }\cdot \nu $ for the push forward of $\nu $
by $\sigma _{\gamma }$. This measure is invariant by $\gamma \phi $. When $%
\gamma =1$, $\sigma _{\gamma }=\sigma _{1}=\sigma $ and $\nu _{1}$ is equal
the measure $\nu _{\Theta _{\mathrm{Mo}}}^{\phi }$ of the last section.

In the proof of Theorem \ref{teoprinc} it will be needed a differentiable
parametrization of the measures $\nu _{\gamma }$. To write it let $\mathbb{E}%
_{\Theta }$ be a flag bundle admitting a section $\sigma $. Consider the
evaluation map $E:\mathcal{G}\rightarrow \Gamma \mathbb{E}_{\Theta }$ given
by 
\begin{equation*}
E\left( \gamma \right) =\gamma \circ \sigma
\end{equation*}%
which is differentiable since it is the action of the gauge group on the
space of sections. The parametrization will be obtained by an application of
the implicit function theorem to $E$ around the identity. For it we must
check that the differential $dE_{1}$ of $E$ at the identity $1\in \mathcal{G}
$ is onto and its kernel is complemented by a closed subspace.

To prove these facts let $f_{\sigma }\in \mathcal{C}_{\mathrm{eq}}(R,\mathbb{%
F}_{\Theta })$ be the equivariant map of the section $\sigma $. It defines
the $K_{\Theta }$-subbundle of $R$, $\mathcal{L}_{\sigma }=f_{\sigma
}^{-1}\{b_{\Theta }\}$ where $b_{\Theta }$ is the origin of $\mathbb{F}%
_{\Theta }$ and $K_{\Theta }$ the isotropy group at $b_{\Theta }$ of the $K$
action on $\mathbb{F}_{\Theta }$. The space of sections $\Gamma \mathbb{E}%
_{\Theta _{\mathrm{Mo}}}$ can be seen as the space of equivariant maps $%
\mathcal{C}_{\mathrm{eq}}(\mathcal{L}_{\sigma },\mathbb{F}_{\Theta })$.
Analogously \ $\mathcal{G}$ is $\mathcal{C}_{\mathrm{eq}}(\mathcal{L}%
_{\sigma },G)$. Their tangent spaces are $T_{\sigma }\mathcal{C}_{\mathrm{eq}%
}(\mathcal{L}_{\sigma },\mathbb{F}_{\Theta })=\mathcal{C}(\mathcal{L}%
_{\sigma },\mathfrak{n}_{\Theta }^{-})$ and $T_{1}\mathcal{C}_{\mathrm{eq}}(%
\mathcal{L}_{\sigma },G)=\mathcal{C}(\mathcal{L}_{\sigma },\mathfrak{g})$.

Now, let $f_{\gamma }:\mathcal{L}_{\sigma }\rightarrow G$ be the equivariant
map of $\gamma \in \mathcal{G}$. Then $f_{\gamma \circ \sigma }\left(
r\right) =f_{\gamma }\left( r\right) \cdot f_{\sigma }\left( r\right) $,
which means that in terms of the equivariant maps $E$ is given by 
\begin{equation*}
E\left( f_{\gamma }\right) =f_{\gamma }\cdot f_{\sigma }.
\end{equation*}%
Its differential is obtained by taking pointwise derivatives. Thus if $Y\in 
\mathcal{C}(\mathcal{L}_{\sigma },\mathfrak{g})$ then 
\begin{equation*}
\left( dE_{1}\left( Y\right) \right) (r)=dp_{1}\cdot Y(r)
\end{equation*}%
where $p:G\rightarrow \mathbb{F}_{\Theta }$ is the projection. We have $%
\mathfrak{g}=\mathfrak{n}_{\Theta }^{-}\oplus \mathfrak{p}_{\Theta }$ where $%
dp_{1}$ is the projection onto the first component. Any map in $\mathcal{C}(%
\mathcal{L}_{\sigma },\mathfrak{n}_{\Theta }^{-})$ has the form $dp_{1}\cdot
Y(r)$, $Y\in \mathcal{C}(\mathcal{L}_{\sigma },\mathfrak{g})$, that is, $%
dE_{1}$ is onto. Its kernel is given by $Y\in \mathcal{C}(\mathcal{L}%
_{\sigma },\mathfrak{g})$ such that $dp_{1}\cdot Y(r)=0$, that is, 
\begin{equation*}
\ker dE_{1}=\{Y\in \mathcal{C}(\mathcal{L}_{\sigma },\mathfrak{g}):\forall
r\in \mathcal{L}_{\sigma },~Y(r)\in \mathfrak{p}_{\Theta }\}.
\end{equation*}%
This subspace is complemented by $\mathcal{C}(\mathcal{L}_{\sigma },%
\mathfrak{n}_{\Theta }^{-})$ which is a closed subspace.

Hence an application of the implicit function theorem to $E$ at $1$ yields
the following parametrization of sections in a neighborhood of $\sigma $.

\begin{proposicao}
\label{propimpfcthm2}Given a section $\sigma \in \Gamma \mathbb{E}_{\Theta }$
there exists a neighborhood $U$ of $\sigma $ in $\Gamma \mathbb{E}_{\Theta }$
and a differentiable map $l:U\rightarrow \mathcal{G}$ such that $\tau
=l\left( \tau \right) \cdot \sigma $ for all $\tau \in U$. Its differential $%
dl_{\sigma }$ at $\sigma $ associates $Y\in T_{\sigma }\mathcal{C}_{\mathrm{%
eq}}(\mathcal{L}_{\sigma },\mathbb{F}_{\Theta })=\mathcal{C}(\mathcal{L}%
_{\sigma },\mathfrak{n}_{\Theta }^{-})$ with itself via the identification $%
T_{1}\mathcal{C}_{\mathrm{eq}}(\mathcal{L}_{\sigma },G)=\mathcal{C}(\mathcal{%
L}_{\sigma },\mathfrak{g})$.
\end{proposicao}

\begin{profe}
By the implicit function theorem there exists a diffeomorphism $\psi
:A\subset \mathcal{G}\rightarrow B\times U\subset \ker dE_{1}\times \Gamma 
\mathbb{E}_{\Theta }$ such that $E\circ \psi ^{-1}$ is the projection onto
the second component. Hence, $l\left( \tau \right) =\psi ^{-1}\left( 0,\tau
\right) $ is the desired map.
\end{profe}

Now it is easy to get the parametrization of the $\gamma \phi $-invariant
sections $\sigma _{\gamma }\in \Gamma \mathbb{E}_{\Theta _{\mathrm{Mo}}}$,
with $\gamma $ running through an identity neighborhood of $\mathcal{G}$.

In fact, if $V\subset \mathcal{G}$ is open with $l\left( U\right) \subset V$
then $L:V\rightarrow \mathcal{G}$,%
\begin{equation*}
L\left( \gamma \right) =l\left( \sigma _{\gamma }\right) =l\circ s\left(
\gamma \right) ,
\end{equation*}%
is well defined, differentiable and satisfies $\sigma _{\gamma }=L\left(
\gamma \right) \cdot \sigma _{1}$.

This parametrization can be carried out to the $\gamma \phi $-invariant
measures $\nu _{\gamma }=\sigma _{\gamma }\cdot \nu $, namely $\nu _{\gamma
}=L\left( \gamma \right) \cdot \nu _{1}$.

Summarizing we get the following ingredient for the proof of Theorem \ref%
{teoprinc}.

\begin{proposicao}
\label{promedinvpertub}There exists an identity neighborhood $V$ of $%
\mathcal{G}$ and a map $L:V\rightarrow \mathcal{G}$ such that $\sigma
_{\gamma }=L\left( \gamma \right) \cdot \sigma _{1}$. If $\nu _{\gamma
}=\sigma _{\gamma }\cdot \nu $ then $\nu _{\gamma }=L\left( \gamma \right)
\cdot \nu _{1}$.

Its differential $dL_{1}=dl_{\sigma }\circ ds_{1}$ assumes values in $%
\mathcal{C}(\mathcal{L}_{\sigma },\mathfrak{n}_{\Theta }^{-})$.
\end{proposicao}

\subsection{Proof of smoothness\label{subsecprfdif}}

To finish the proof of Theorem \ref{teoprinc} it remains to check the
differentiability of the map 
\begin{equation*}
\gamma \in \mathcal{G}\longmapsto \int \mathsf{a}_{\omega }^{\gamma \phi
}\left( 1,\xi \right) \nu _{\gamma }\left( d\xi \right)
\end{equation*}%
(see Proposition \ref{propintflag}). By Proposition \ref{promedinvpertub} we
have $\nu _{\gamma }=L\left( \gamma \right) \cdot \nu _{1}$ hence the above
integral over $\nu _{\gamma }$ equals to 
\begin{equation*}
\int \mathsf{a}_{\omega }^{\gamma \phi }\left( 1,L\left( \gamma \right)
\left( \xi \right) \right) \nu _{1}\left( d\xi \right) ,
\end{equation*}%
where $L\left( \gamma \right) \left( \xi \right) $ stands for the action of $%
L\left( \gamma \right) \in \mathcal{G}$ on $\xi \in \mathbb{E}_{\Theta _{%
\mathrm{Mo}}\left( \phi \right) }$. Now, consider the map $\Upsilon ^{\omega
}:\mathcal{G}\rightarrow \mathcal{C}\left( \mathbb{E}_{\Theta _{\mathrm{Mo}%
}\left( \phi \right) },\mathbb{R}\right) $ defined by%
\begin{equation*}
\Upsilon ^{\omega }\left( \gamma \right) \left( \cdot \right) =\mathsf{a}%
_{\omega }^{\gamma \phi }\left( 1,L\left( \gamma \right) \left( \cdot
\right) \right) .
\end{equation*}%
This map is differentiable. In fact by (\ref{fora}) we have 
\begin{equation*}
\mathsf{a}_{\omega }^{\gamma \phi L\left( \gamma \right) }\left( 1,\xi
\right) =\mathsf{a}_{\omega }^{\gamma \phi }\left( 1,L\left( \gamma \right)
\left( \xi \right) \right) +\mathsf{a}_{\omega }^{L\left( \gamma \right)
}\left( 1,\xi \right) ,
\end{equation*}%
hence 
\begin{equation*}
\Upsilon ^{\omega }\left( \gamma \right) \left( \xi \right) =\mathsf{a}%
_{\omega }^{\gamma \phi L\left( \gamma \right) }\left( 1,\xi \right) -%
\mathsf{a}_{\omega }^{L\left( \gamma \right) }\left( 1,\xi \right) .
\end{equation*}%
By Proposition \ref{propdifcocycle} both terms in the second hand side are
differentiable as functions of $\gamma $. The second term is $F_{\omega }^{%
\mathrm{id}}\circ L$, while the first term is the composition of $F_{\omega
}^{\phi }$ with the map 
\begin{equation}
p:\gamma \mapsto \gamma C_{\phi }\left( L\left( \gamma \right) \right) ,
\label{foraplicapeconj}
\end{equation}%
where $C_{\phi }\left( \eta \right) =\phi \eta \phi ^{-1}$. This map is
differentiable as well because $\mathcal{G}$ is a Banach Lie group and
conjugation by $\phi $ is differentiable.

Finally, 
\begin{equation*}
\int \mathsf{a}_{\omega }^{\gamma \phi }\left( 1,\xi \right) \nu _{\gamma
}\left( d\xi \right) =\int \Upsilon ^{\omega }\left( \gamma \right) \left(
\xi \right) \nu _{1}\left( d\xi \right)
\end{equation*}%
is the composition of a continuous linear map on $\mathcal{C}\left( \mathbb{E%
}_{\Theta _{\mathrm{Mo}}},\mathbb{R}\right) $ (the measure $\nu _{1}$) with
the differentiable map $\Upsilon ^{\omega }$. This proves the smoothness
statement of Theorem \ref{teoprinc}.

\subsection{Differential}

The differentiability of the partial Lyapunov spectrum 
\begin{equation*}
\Lambda _{\omega }\left( \gamma \right) =\int_{X}\omega \left( H_{\gamma
\phi }^{+}\left( x\right) \right) \nu \left( dx\right)
\end{equation*}%
was proved by showing that this map is a composition of several
differentiable maps whose differentials we have recollected above. Now we
can compose these differentials to get the differential of $\Lambda _{\omega
}$ at the identity $1\in \mathcal{G}$.

\begin{proposicao}
The differential of $\Lambda _{\omega }\left( \gamma \right) =\int_{X}\omega
\left( H_{\gamma \phi }^{+}\left( x\right) \right) \nu \left( dx\right) $ at 
$\gamma =1$ is 
\begin{equation}
d\left( \Lambda _{\omega }\right) _{1}\left( Y\right) =\int_{\mathbb{E}%
_{\Theta }}\omega \circ \mathbf{a}\left( f_{Y}\left( \phi ^{R}\left(
r\right) \right) \right) \left( \sigma \cdot \nu \right) \left( d\xi \right)
\qquad \xi =r\cdot b_{\Theta }  \label{fordiferenexpo}
\end{equation}%
where $Y\in T_{1}\mathcal{G}$ is a section of the bundle $Q\times _{\mathrm{%
Ad}}\mathfrak{g}$, $f_{Y}:R\rightarrow \mathfrak{g}$ its equivariant map $Y$
and $\mathbf{a}:\mathfrak{g}\rightarrow \mathfrak{a}$ is the projection onto
the $\mathfrak{a}$-component in the Iwasawa decomposition $\mathfrak{g}=%
\mathfrak{k}\oplus \mathfrak{a}\oplus \mathfrak{n}$.
\end{proposicao}

\begin{profe}
As settled above $\Lambda _{\omega }\left( \gamma \right) $ is the integral
of the map $\Upsilon ^{\omega }:\mathcal{G}\rightarrow \mathcal{C}\left( 
\mathbb{E}_{\Theta _{\mathrm{Mo}}\left( \phi \right) },\mathbb{R}\right) $
given by 
\begin{equation*}
\Upsilon ^{\omega }\left( \gamma \right) \left( \cdot \right) =\mathsf{a}%
_{\omega }^{\gamma \phi L\left( \gamma \right) }\left( 1,\cdot \right) -%
\mathsf{a}_{\omega }^{L\left( \gamma \right) }\left( 1,\cdot \right)
\end{equation*}%
with $\Upsilon ^{\omega }=F_{\omega }^{\phi }\circ p-F_{\omega }^{\mathrm{id}%
}\circ L$ where $F_{\omega }^{\phi }:\mathcal{G}\rightarrow \mathcal{C}%
\left( \mathbb{E}_{\Theta },\mathfrak{a}_{\Theta }\right) $ were defined in (%
\ref{fordefFtetafi}) and $p:\mathcal{G}\rightarrow \mathcal{G}$ is $p\left(
\gamma \right) =\gamma \cdot \left( C_{\phi }\circ L\right) \left( \gamma
\right) $. Hence the differential of $\Lambda _{\omega }$ is the integral of
the differential of $F_{\omega }^{\phi }\circ p-F_{\omega }^{\mathrm{id}%
}\circ L$.

Now by Proposition \ref{propdifcocycle} we have 
\begin{equation}
d\left( F_{\omega }^{\phi }\right) _{1}\left( Y\right) \left( r\cdot
b_{\Theta }\right) =\omega \circ \mathbf{a}\left( f_{Y}\circ \phi ^{R}\left(
r\right) \right) .  \label{fordifdenovo}
\end{equation}%
Since $dL_{1}$ assumes values in $\mathcal{C}(\mathcal{L}_{\sigma },%
\mathfrak{n}_{\Theta }^{-})$, it follows that the differential at $1$ of $%
F_{\omega }^{\mathrm{id}}\circ L$ is zero. As to the term $F_{\omega }^{\phi
}\circ p$ we have first 
\begin{equation*}
dp_{1}\left( Y\right) =Y+d\left( C_{\phi }\circ L\right) _{1}\left( Y\right)
.
\end{equation*}%
The $\mathfrak{a}$-component in the Iwasawa decomposition of the second term \linebreak
$d\left( C_{\phi }\circ L\right) _{1}\left( Y\right) $ annihilates, since $%
dL_{1}$ assumes values in $\mathcal{C}(\mathcal{L}_{\sigma },\mathfrak{n}%
_{\Theta }^{-})$ and $d\left( C_{\phi }\right) _{1}$ is just composition of
the equivariant map with $\phi ^{-1}$. Hence, for any $Y\in \mathrm{Lie}%
\left( \mathcal{G}\right) $ the equivariant map of $d\left( C_{\phi }\circ
L\right) _{1}\left( Y\right) $ takes values in $\mathfrak{n}_{\Theta }^{-}$
so that its $\mathfrak{a}$-component is zero.

We are left with $d\left( F_{\omega }^{\phi }\right) _{1}\left( Y\right) $
which is $\omega \circ \mathbf{a}\left( f_{Y}\circ \phi ^{R}\right) $ by (%
\ref{fordifdenovo}).
\end{profe}

To use below in examples we rewrite formula (\ref{fordiferenexpo}) in the
case of a trivial bundle $X\times Q$ where the integral can be taken with
respect to $\nu $ on the base $X$. Here $\phi ^{n}\left( x,g\right) =\left(
\tau ^{n}\left( x\right) ,\rho \left( n,x\right) g\right) $ is defined by a
cocycle $\rho \left( n,x\right) $ generated by the map $\rho \left(
1,x\right) \in G$. We take $R=X\times K$, so that $\phi ^{R}\left(
x,k\right) =\left( \tau \left( x\right) ,u\left( x,k\right) \right) $ where $%
\rho \left( 1,x\right) k=u\left( x,k\right) hn$ is the Iwasawa decomposition.

The element $Y\in \mathrm{Lie}\left( \mathcal{G}\right) $ is a map $%
Y:X\rightarrow \mathfrak{g}$ that equals the restriction of $f_{Y}$ to $%
X\times \{1\}$ so that $f_{Y}\left( x,k\right) =\mathrm{Ad}\left(
k^{-1}\right) Y\left( x\right) $.

Let $\sigma :X\rightarrow \mathbb{F}_{\Theta }$ be an invariant section and
choose $k_{x}\in K$ such that $\sigma \left( x\right) =k_{x}\left( b_{\Theta
}\right) $. Then 
\begin{equation*}
f_{Y}\left( \phi ^{R}\left( x,k_{x}\right) \right) =f_{Y}\left( \tau \left(
x\right) ,u\left( x,k_{x}\right) \right) =\mathrm{Ad}\left( u\left(
x,k_{x}\right) ^{-1}\right) Y\left( \tau \left( x\right) \right) .
\end{equation*}%
It follows by (\ref{fordiferenexpo}) that 
\begin{equation}
d\left( \Lambda _{n\omega }\right) _{1}\left( Y\right) =\int_{X}\omega \circ 
\mathbf{a}\left( \mathrm{Ad}\left( u\left( x,k_{x}\right) ^{-1}\right)
Y\left( \tau \left( x\right) \right) \right) \nu \left( dx\right) .
\label{fordiferentrivial}
\end{equation}%
The choice of $k_{x}$ is immaterial.

\section{A general differentiable perturbation\label{secgener}}

So far we proved in Theorem \ref{teoprinc} differentiability of the Lyapunov
spectra with respect to perturbations by the gauge group $\mathcal{G}$. In
this section we abstract the properties of $\mathcal{G}$ required in the
proof of differentiability and state a more general set up where the proof
of Theorem \ref{teoprinc} still applies.

Thus we let $\mathcal{H}\subset \mathrm{Aut}\left( Q\right) $ be a
Banach-Lie group of automorphisms of $Q$ and take $\phi \in \mathrm{Aut}%
\left( Q\right) $. Then the proof of Theorem \ref{teoprinc} holds for $%
\mathcal{H}$ in place of $\mathcal{G}$ if the following conditions are
satisfied.

\begin{itemize}
\item[\textbf{H1}] The maps induced on the base space $X$ by $\phi $ and any 
$\gamma \in \mathcal{H}$ preserves a probability measure $\nu $.

\item[\textbf{H2}] The maps $\mathcal{H\phi }\times \Gamma \mathbb{E}%
_{\Theta _{\mathrm{Mo}}\left( \phi \right) }\rightarrow \Gamma \mathbb{E}%
_{\Theta _{\mathrm{Mo}}\left( \phi \right) }$ and $\mathcal{H\phi }\times
\Gamma \mathbb{E}_{\Theta _{\mathrm{Mo}}^{\ast }\left( \phi \right)
}\rightarrow \Gamma \mathbb{E}_{\Theta _{\mathrm{Mo}}^{\ast }\left( \phi
\right) }$ given by the actions of $\mathcal{H\phi }$ on the spaces of
sections $\Gamma \mathbb{E}_{\Theta _{\mathrm{Mo}}\left( \phi \right) }$ and 
$\Gamma \mathbb{E}_{\Theta _{\mathrm{Mo}}^{\ast }\left( \phi \right) }$ are
smooth.

\item[\textbf{H3}] $\mathcal{H}$ is normalized by $\phi $ and the
conjugation $C_{\phi }:\mathcal{H}\rightarrow \mathcal{H}$, $C_{\phi }\left(
\eta \right) =\phi \eta \phi ^{-1}$ is smooth.

\item[\textbf{H4}] If $\omega $ belongs to the linear span of $\Omega
\setminus \Omega _{\Theta _{\mathrm{Mo}}\left( \phi \right) }$ then the map 
\begin{equation*}
F_{\omega }^{\phi }:\mathcal{H}\rightarrow \mathcal{C}\left( \mathbb{E}%
_{\Theta },\mathbb{R}\right) \qquad F_{\omega }^{\phi }\left( \gamma \right)
=\mathsf{a}_{\omega }^{\gamma \phi }\left( 1,\cdot \right)
\end{equation*}%
\qquad is smooth.

\item[\textbf{H5}] Let $\sigma _{\mathrm{Mo}}\in \Gamma \mathbb{E}_{\Theta _{%
\mathrm{Mo}}\left( \phi \right) }$ be the $\phi $-invariant section whose
image is the attractor Morse component of $\phi $ on $\mathbb{E}_{\Theta _{%
\mathrm{Mo}}\left( \phi \right) }$. Then $\mathcal{H}$ is locally transitive
at $\sigma _{\mathrm{Mo}}$ in the following sense: Define the map $E:%
\mathcal{H}\rightarrow \Gamma \mathbb{E}_{\Theta _{\mathrm{Mo}}\left( \phi
\right) }$ by $E\left( \gamma \right) =\gamma \left( \sigma _{\mathrm{Mo}%
}\right) $ which smooth by (H2). Then it is required the differential $%
dE_{1} $ is onto and $\ker dE_{1}$ is complementable in $T_{1}\mathcal{H}$.
\end{itemize}

Clearly (H3) is automatic if $\mathcal{H}$ contains $\phi $.

To see that Theorem \ref{teoprinc} can be proved under these conditions,
note first that by (H2) we can apply the Implicit Function Theorem as in
Proposition \ref{propptofix} to get $\gamma \phi $-invariant sections with $%
\gamma $ close to the identity of $\mathcal{H}$. The proof that the
invariant sections on $\Gamma \mathbb{E}_{\Theta _{\mathrm{Mo}}\left( \phi
\right) }$ give the attractor Morse components (Theorem \ref{teosecMorse})
rely only on properties of the Morse decompositions.

Now smoothness of the cocycles, as in Section \ref{secociclan}, is
axiomatized in (H4) while condition (H5) allows a new application of the
Implicit Function Theorem as in Proposition \ref{propimpfcthm2} to get a
smooth parametrization of the invariant measures on $\mathbb{E}_{\Theta _{%
\mathrm{Mo}}\left( \phi \right) }$.

Finally the assumption that $\mathcal{H}$ is a Banach-Lie group and (H3) are
required at the end of the proof as in Subsection \ref{subsecprfdif} to
ensure that the map $p:\gamma \mapsto \gamma C_{\phi }\left( L\left( \gamma
\right) \right) $ defined in (\ref{foraplicapeconj}) is smooth.

\section{Cocycles and semigroups\label{secsemi}}

In this section we let $\rho :\mathbb{Z}\times X\rightarrow G$ be a cocycle
over the homeomorphism $\tau :X\rightarrow X$ leaving invariant the measure $%
\nu $. It defines the skew-product flow $\phi ^{n}\left( x,g\right) =\left(
\tau \left( x\right) ,\rho \left( n,x\right) g\right) $ on $X\times G$. In
case $\rho $ takes values in a semigroup $S\subset G$ then the flag type $%
\Theta _{\mathrm{Mo}}\left( \phi \right) $ of $\phi $ can be obtained in
case we know in advance the flag type of $S$. This is the set up of several
examples in the literature (see Barreira-Pesin \cite{barrpe}). Also in
Ruelle \cite{ruanl} the flow is given this way with $G=\mathrm{Gl}\left( d,%
\mathbb{R}\right) $ and $S\subset G$ the semigroup of linear maps leaving
invariant a cone $W\subset \mathbb{R}^{d}$.

We refer to \cite{sminv}, \cite{SM}, \cite{SMT}, \cite{smmax} and references
therein to results about semigroups in semi-simple (and reductive) Lie
groups and their actions on flag manifolds. In \cite{smbflow} and \cite{msm}
these semigroups were related to flows on principal bundles via the
shadowing semigroups. The main result of these papers is that the Morse
components on the flag bundles are given by shriking the control sets of the
shadowing semigroups.

What concern us here is the flag type of $S\subset G$ (when $\mathrm{int}%
S\neq \emptyset $) which is the largest flag manifold $\mathbb{F}_{\Theta
\left( S\right) }$ of $G$ such that the unique invariant control set of $S$
in $\mathbb{F}_{\Theta \left( S\right) }$ is contractible (in the sense that
is contained in an open Bruhat cell).

If the cocycle $\rho $ is contained in $\mathrm{int}S$ then we prove below
that the flag type $\Theta _{\mathrm{Mo}}\left( \phi \right) $ contains the
flag type $\Theta \left( S\right) $ of $S$, so that the flag manifolds
themselves $\mathbb{F}_{\Theta _{\mathrm{Mo}}\left( \phi \right) }$ and $%
\mathbb{F}_{\Theta \left( S\right) }$ are related by an equivariant
fibration $\mathbb{F}_{\Theta _{\mathrm{Mo}}\left( \phi \right) }\rightarrow 
\mathbb{F}_{\Theta \left( S\right) }$, that is, $\mathbb{F}_{\Theta _{%
\mathrm{Mo}}\left( \phi \right) }$ is larger than $\mathbb{F}_{\Theta \left(
S\right) }$.

For several semigroups $S$ the flag type $\Theta \left( S\right) $ is known
so that from the inclusion $\Theta _{\mathrm{Mo}}\left( \phi \right) \subset
\Theta \left( S\right) $ we get smoothness of the part of the $\mathfrak{a}$%
-Lyapunov spectra given by $\omega \in \mathrm{span}\left( \Omega \setminus
\Omega _{\Theta \left( S\right) }\right) $. The inclusion $\Theta _{\mathrm{%
Mo}}\left( \phi \right) \subset \Theta \left( S\right) $ also implies the
positivity of some Lyapunov exponents by the results of \cite{smsec}.

From now on let $\rho :\mathbb{Z}\times X\rightarrow G$ be a continuous
cocycle such that $\rho \left( 1,x\right) \in \mathrm{int}S$ (and hence $%
\rho \left( n,x\right) \in \mathrm{int}S$) for all $x\in X$ where $S$ is a
subsemigroup of $G$. If $\mathbb{F}_{\Theta \left( S\right) }$ is the flag
type of $S$ then we write $C_{\Theta \left( S\right) }$ (or simply $C$) for
the unique invariant control set of $S$ on $\mathbb{F}_{\Theta \left(
S\right) }$ (see \cite{sminv}). It is compact and has nonempty interior.

To prove the $\Theta _{\mathrm{Mo}}\left( \phi \right) \subset \Theta \left(
S\right) $ we need the following lemma where we use the metric space
notation $d\left( K_{1},K_{2}\right) =\inf \{d\left( x_{1},x_{2}\right)
:x_{i}\in K_{i}\}$ and $d\left( x,K_{1}\right) =\inf \{d\left( x,y\right)
:y\in K_{1}\}$ if $K_{1}$ and $K_{2}$ are compact and $x$ is a point.

\begin{lema}
There exists $\varepsilon >0$ such that $d\left( g_{x}C,\mathrm{cl}%
C^{c}\right) >\varepsilon $ for all $x\in X$ where $g_{x}=\rho \left(
1,x\right) $ and $C=C_{\Theta \left( S\right) }$.
\end{lema}

\begin{profe}
Since $g_{x}\in \mathrm{int}S$, we have $g_{x}C\subset \mathrm{int}C$ and
hence there exists $\varepsilon _{x}>0$ such that $d\left( g_{x}C,\mathrm{cl}%
C^{c}\right) >\varepsilon _{x}$. By continuity with respect to the
compact-open topology there exists an open set $V_{x}\subset G$ with $%
g_{x}\in V_{x}$ such that $d\left( gC,\mathrm{cl}C^{c}\right) >\varepsilon
_{x}$ if $g\in V_{x}$. Hence by compactness of $X$ we get $\varepsilon >0$
such that $d\left( g_{x}C,\mathrm{cl}C^{c}\right) >\varepsilon $ for all $%
x\in X$.
\end{profe}

If $\varepsilon >0$ is as in this lemma then any $\left( \varepsilon
/2,T\right) $-chain $\left( x_{0},\xi _{0}\right) $, \ldots , $\left(
x_{n},\xi _{n}\right) $ starting at $\left( x_{0},\xi _{0}\right) \in
X\times C$ is contained in $X\times C$. In fact, $\left( x_{1},\xi
_{1}\right) $ is $\varepsilon /2$ close to $\phi ^{t}\left( x_{0},\xi
_{0}\right) $. But $\phi ^{t-1}\left( x_{0},\xi _{0}\right) \in X\times C$
so that $\phi ^{t}\left( x_{0},\xi _{0}\right) $ belongs to $X\times \left(
g_{\tau ^{t-1}\left( x_{0}\right) }C\right) $. Hence $\xi _{1}$ is $%
\varepsilon /2$ close to $g_{\tau ^{t-1}\left( x_{0}\right) }C$, implying
that $\xi _{1}\in C$. By induction it follows that the whole chain is
contained in $X\times C$. \ 

As a consequence we have that the $\omega $-limit $\omega \left( x,\xi
\right) $ is contained in $X\times C$ if $\left( x,\xi \right) \in X\times C$%
.

Now, let $\mathcal{M}_{\Theta \left( S\right) }^{+}$ be the attractor Morse
component on $X\times \mathbb{F}_{\Theta \left( S\right) }$. Then for a
dense subset $D\subset X\times \mathbb{F}_{\Theta \left( S\right) }$ we have 
$\omega \left( x,\xi \right) \subset \mathcal{M}_{\Theta \left( S\right)
}^{+}$ if $\left( x,\xi \right) \in D$. Since $X\times C$ has nonempty
interior we have $\left( X\times C\right) \cap D\neq \emptyset $ so that $%
\omega \left( x,\xi \right) \subset \left( X\times C\right) \cap \mathcal{M}%
_{\Theta \left( S\right) }^{+}$ if $\left( x,\xi \right) \in \left( X\times
C\right) \cap D$. Thus we conclude that $\left( X\times C\right) \cap 
\mathcal{M}_{\Theta \left( S\right) }^{+}\neq \emptyset $. Finally we get $%
\mathcal{M}_{\Theta \left( S\right) }^{+}\subset X\times C$ since $\mathcal{M%
}_{\Theta \left( S\right) }^{+}$ is chain transitive and $X\times C$ is
chain invariant.

\begin{proposicao}
\label{propflagcontido}If $\rho \left( 1,x\right) \in \mathrm{int}S$, $x\in
X $, then the attractor Morse component $\mathcal{M}_{\Theta \left( S\right)
}^{+}$ is contained in $X\times C$ and $\Theta _{\mathrm{Mo}}\left( \phi
\right) \subset \Theta \left( S\right) $.
\end{proposicao}

\begin{profe}
The first statement was proved above. The second statement is a consequence
of the first one since the inclusion $\mathcal{M}_{\Theta \left( S\right)
}^{+}\subset X\times C$ implies that for any $x\in X$ the fiber $\mathcal{M}%
_{\Theta \left( S\right) }^{+}\left( x\right) =\mathcal{M}_{\Theta \left(
S\right) }^{+}\cap \left( \{x\}\times \mathbb{F}_{\Theta \left( S\right)
}\right) $ is contained in an open Bruhat cell. This forces $\mathcal{M}%
_{\Theta \left( S\right) }^{+}\left( x\right) $ to be a point since this is
the case where the algebraic manifold $\mathcal{M}_{\Theta \left( S\right)
}^{+}\left( x\right) $ is contained in an open Bruhat cell.
\end{profe}

As mentioned above the inclusion $\Theta _{\mathrm{Mo}}\left( \phi \right)
\subset \Theta \left( S\right) $ implies differentiability and positivity of
certain linear maps of the $\mathfrak{a}$-Lyapunov exponents.

\begin{corolario}
\label{corsemigroupres}Suppose that $\rho \left( 1,x\right) \in \mathrm{int}%
S $ and let $\Theta \left( S\right) $ be the flag type of $S$. Then

\begin{enumerate}
\item $\Lambda _{\omega }\left( \gamma \right) =\int_{X}\omega \left(
H_{\gamma \phi }^{+}\left( x\right) \right) \nu \left( dx\right) $ is smooth
if $\omega \in \mathrm{span}\left( \Omega \setminus \Omega _{\Theta \left(
S\right) }\right) $.

\item If $\alpha \notin \Theta \left( S\right) $ then $\alpha \left( H_{\phi
}^{+}\left( x\right) \right) >0$ for $x\in X$ where $H_{\phi }^{+}\left(
x\right) $ is defined.
\end{enumerate}
\end{corolario}

\begin{profe}
The first statement is a consequence of our Theorem \ref{teoprinc}, while
the second one is due to the description of the Morse spectra made in \cite%
{smsec} as a convex set contained in 
\begin{equation*}
\bigcap\limits_{\alpha \notin \Theta \left( S\right) }\{H\in \mathfrak{a}%
:\alpha \left( H\right) >0\}.
\end{equation*}
\end{profe}

In the sequel we give some examples of semigroups together with their flag
types and interpret the above corollary on the light of these flags.

\subsection{Invariant cones}

Let $W\subset \mathbb{R}^{d}$ be a pointed (i.e., $W\cap -W=\{0\}$) and
generating (i.e., $\mathrm{int}W\neq \emptyset $) cone. The semigroup%
\begin{equation*}
S_{W}=\{g\in \mathrm{Sl}\left( d,\mathbb{R}\right) :gW\subset W\cup -W\}
\end{equation*}%
has nonempty interior and its flag type is the projective space $\mathbb{P}%
^{d-1}$ (see \cite{SMT}, Example 5.3). For the group $\mathrm{Sl}\left( d,%
\mathbb{R}\right) $ and its Lie algebra $\mathfrak{sl}\left( d,\mathbb{R}%
\right) $ we have $\mathfrak{a}\subset \mathfrak{sl}\left( d,\mathbb{R}%
\right) $ is the abelian algebra of diagonal matrices. The roots are $\alpha
_{ij}=\lambda _{i}-\lambda _{j}$, $i\neq j$, where 
\begin{equation*}
\lambda _{i}\left( \mathrm{diag}\{a_{1},\ldots ,a_{d}\}\right) =a_{i}.
\end{equation*}%
A natural choice of positive roots is $\Pi ^{+}=\{\alpha _{ij}:i<j\}$ that
have the set of simple roots 
\begin{equation*}
\Sigma =\{\alpha _{12},\alpha _{23},\ldots ,\alpha _{d-1,d}\}.
\end{equation*}%
The associated fundamental weights are 
\begin{equation*}
\Omega =\{\lambda _{1},\lambda _{1}+\lambda _{2},\ldots ,\lambda _{1}+\cdots
+\lambda _{d-1}\}
\end{equation*}

The subset $\Theta =\Theta \left( S_{W}\right) \subset \Sigma $
corresponding to $\mathbb{P}^{d-1}$ is 
\begin{equation*}
\Theta \left( S_{W}\right) =\{\alpha _{23},\ldots ,\alpha _{d-1,d}\}
\end{equation*}%
so that $\Omega \setminus \Omega _{\Theta }=\{\lambda _{1}\}$.

Hence if $\rho $ is a cocycle with values in $\mathrm{int}S_{W}$ then by (2)
of Corollary \ref{corsemigroupres}, we have $\alpha _{12}\left( H^{+}\left(
x\right) \right) >0$ if $H^{+}\left( x\right) $ is defined implying that the
largest Lyapunov exponent $\lambda _{1}\left( H^{+}\left( x\right) \right) $
has multiplicity one. By Theorem \ref{teoprinc}, $\rho $ is a differentiable
point of the largest exponent 
\begin{equation*}
\int_{X}\lambda _{1}\left( H^{+}\left( x\right) \right) \nu \left( dx\right)
.
\end{equation*}%
Thus we recover the result of Ruelle \cite{ruanl} on invariant cones.

To write the differential of the largest Lyapunov exponent at $\rho $, let $%
\sigma :X\rightarrow \mathbb{P}^{d-1}$ be the invariant section giving the
attractor Morse component. Its image is contained in $\mathrm{int}W$ so
there is $\overline{\sigma }:X\rightarrow S^{d-1}$ such that $\sigma \left(
x\right) =\mathrm{span}\overline{\sigma }\left( x\right) $. We have 
\begin{equation*}
\phi ^{R}\left( x,\overline{\sigma }\left( x\right) \right) =\frac{\rho
\left( 1,x\right) \overline{\sigma }\left( x\right) }{\left\Vert \rho \left(
1,x\right) \overline{\sigma }\left( x\right) \right\Vert }.
\end{equation*}

In the notation of formula (\ref{fordiferentrivial}) we have $K=\mathrm{SO}%
\left( d\right) $ and $u\left( x,k_{x}\right) e_{1}=\frac{\rho \left(
1,x\right) \overline{\sigma }\left( x\right) }{\left\Vert \rho \left(
1,x\right) \overline{\sigma }\left( x\right) \right\Vert }$ where $e_{1}$ is
the first basis vector so $\left[ e_{1}\right] $ is the origin of $\mathbb{P}%
^{d-1}$. Hence if $Y:X\rightarrow \mathfrak{g}$ is an element of $\mathrm{Lie%
}\left( \mathcal{G}\right) $ then 
\begin{equation*}
\lambda _{1}\circ \mathbf{a}\left( \mathrm{Ad}\left( u\left( x,k_{x}\right)
^{-1}\right) Y\left( x\right) \right) =\frac{\langle Y\left( x\right) \rho
\left( 1,x\right) \overline{\sigma }\left( x\right) ,\rho \left( 1,x\right) 
\overline{\sigma }\left( x\right) \rangle }{\left\Vert \rho \left(
1,x\right) \overline{\sigma }\left( x\right) \right\Vert ^{2}}.
\end{equation*}%
Therefore, 
\begin{equation*}
d\left( \Lambda _{\omega }\right) _{1}\left( Y\right) =\int_{X}\frac{\langle
Y\left( x\right) \rho \left( 1,x\right) \overline{\sigma }\left( x\right)
,\rho \left( 1,x\right) \overline{\sigma }\left( x\right) \rangle }{%
\left\Vert \rho \left( 1,x\right) \overline{\sigma }\left( x\right)
\right\Vert ^{2}}\nu \left( dx\right) .
\end{equation*}

These facts hold for $W=\mathbb{R}_{+}^{d}$ when $S_{W}$ is the semigroup of
matrices with nonnegative entries while $g\in \mathrm{int}S_{W}$ has
strictly positive entries.

\subsection{Totally positive matrices}

Let $\mathcal{T}\subset \mathrm{Sl}\left( d,\mathbb{R}\right) $ be the
semigroup of totally positive matrices, that is, $g\in \mathcal{T}$ if and
only if all its minors as $\geq 0$. The elements in $\mathrm{int}\mathcal{T}$
have strictly positive minors. The flag type of $\mathcal{T}$ is $\Theta
\left( \mathcal{T}\right) =\emptyset $, which means that $\mathbb{F}_{\Theta
\left( \mathcal{T}\right) }$ is the manifold of complete flags of subspaces
of $\mathbb{R}^{d}$ (see \cite{smmax}, Section 6.3.1). Hence if $\rho $ is a
cocycle with values in $\mathrm{int}\mathcal{T}$ then its Lyapunov spectra $%
H_{\rho }^{+}\left( x\right) =\mathrm{diag}\{\chi _{1}\left( x\right)
,\ldots ,\chi _{d}\left( x\right) \}$ is simple $\chi _{1}\left( x\right)
>\cdots >\chi _{d}\left( x\right) $ and $\rho $ is a differentiable point of 
$\int_{X}H_{\rho }^{+}\left( x\right) \nu \left( dx\right) $.

The same holds for cocycles with values in the interior of the semigroup $%
\mathcal{S}$ of sign regular matrices ($g\in \mathcal{S}$ if all its minors
have the same sign), since this semigroup also has flag type $\Theta \left( 
\mathcal{S}\right) =\emptyset $.

As a variation of this example take a subset of integers $\mathbf{k}%
=\{k_{1},\ldots ,k_{s}\}$ with $1\leq k_{1}<\cdots <k_{s}<d$, and define the
semigroup $\mathcal{T}_{\mathbf{k}}$ of matrices whose minors of order $%
k_{1},\ldots ,k_{s}$ are nonnegative. Its flag type is the manifold of flags
of subspaces $\left( V_{1}\subset \cdots \subset V_{s}\right) $ with $\dim
V_{i}=k_{i}$ (see \cite{smmax}, Section 6.3.1). The subset of simple roots $%
\Theta \left( \mathcal{T}_{\mathbf{k}}\right) $ associated to this flag
manifold is 
\begin{equation*}
\Theta \left( \mathcal{T}_{\mathbf{k}}\right) =\Sigma \setminus \{\alpha
_{i,i+1}:i\in \mathbf{k}\}
\end{equation*}%
where $\Sigma =\{\alpha _{12},\alpha _{23},\ldots ,\alpha _{d-1,d}\}$ were
defined above. The fundamental weights corresponding to the complement $%
\{\alpha _{i,i+1}:i\in \mathbf{k}\}$ of $\Theta \left( \mathcal{T}_{\mathbf{k%
}}\right) $ are $\omega _{k_{i}}=\lambda _{1}+\cdots +\lambda _{k_{i}}$ with 
$k_{i}\in \mathbf{k}$. These are the weights in the statement of Theorem \ref%
{teoprinc}. Hence if $\rho $ is a cocycle with values in $\mathrm{int}%
\mathcal{T}_{\mathbf{k}}$ then $\rho $ is a differentiable point of 
\begin{equation*}
\int_{X}\left( \chi _{1}\left( x\right) +\cdots +\chi _{k_{i}}\left(
x\right) \right) \nu \left( dx\right) \qquad k_{i}\in \mathbf{k}.
\end{equation*}

\subsection{Compression semigroups and $\mathcal{B}$-convex sets}

A way of producing subsemigroups of $G$ is by taking a subset $C$ of a flag
manifold $\mathbb{F}_{\Theta }$ and define the compression semigroup 
\begin{equation*}
S_{C}=\{g\in G:gC\subset C\}.
\end{equation*}
Conditions on $C$ to have $\mathrm{int}S_{C}\neq \emptyset $ were obtained
in \cite{smmax}. Following \cite{smmax} we say that a subset $C\subset 
\mathbb{F}_{\Theta }$ is admissible if $C$ is contained in an open Bruhat
cell of $\mathbb{F}_{\Theta }$. Then we have the following fact proved in 
\cite{smmax}, Proposition 4.2.

\begin{proposicao}
\label{propcompression}If $C\subset \mathbb{F}_{\Theta }$ is admissible and $%
C=\mathrm{cl}\left( \mathrm{int}C\right) $ then $\mathrm{int}S_{C}=\{g\in
G:gC\subset \mathrm{int}C\}\neq \emptyset $, $C$ is the invariant control
set of $S_{C}$ and $\mathbb{F}_{\Theta }$ is the flag type of $S_{C}$.
\end{proposicao}

Combining this fact with Proposition \ref{propflagcontido} above we have the
following upper bound for the flag type of a cocycle.

\begin{proposicao}
\label{propcompressionro}Let $C\subset \mathbb{F}_{\Theta }$ be an
admissible subset with $C=\mathrm{cl}\left( \mathrm{int}C\right) $ and
suppose that the cocycle $\rho $ in $G$ satisfies $\rho \left( 1,x\right)
C\subset \mathrm{int}C$ for all $x\in X$ then $\Theta \left( \phi \right)
\subset \Theta $, and the conclusions of Corollary \ref{corsemigroupres}
hold with $\Theta $ in place of $\Theta \left( S\right) $.
\end{proposicao}

As an example of a flag manifold we mention the real projective space $%
\mathbb{R}P^{d-1}$ as a flag manifold of $\mathrm{Sl}\left( d,\mathbb{R}%
\right) $. A Bruhat cell is the complement of a hyperplane $\left[ V\right]
=\{\left[ v\right] \in \mathbb{R}P^{d-1}:v\in V\}$ where $V\subset \mathbb{R}%
^{d}$ has $\dim V=d-1$. If $W\subset \mathbb{R}^{d}$ is a pointed and
generating convex cone then $C=\left[ W\right] =\{\left[ v\right] \in 
\mathbb{R}P^{d-1}:v\in W\}$ is admissible and satisfies $C=\mathrm{cl}\left( 
\mathrm{int}C\right) $. In this case $S_{W}=S_{C}$ hence the semigroups
leaving invariant cones treated above are special cases of this set up.

When $d=2$ a hyperplane in $\mathbb{R}P^{1}$ reduces to a point so that a
subset $C=\mathrm{cl}\left( \mathrm{int}C\right) $ is admissible if and only
if it is proper. Hence the above statements yield the following fact.

\begin{proposicao}
Let $\rho $ be a cocycle in $\mathrm{Sl}\left( 2,\mathbb{R}\right) $ such
that $\rho \left( 1,x\right) C\subset \mathrm{int}C$ (all $x\in X$) for a
subset $C$ of $\mathbb{R}P^{1}$ with $C=$. Then,

\begin{enumerate}
\item there exists $c>0$ such that for almost all $x\in X$ the Lyapunov
exponents are $\xi _{\rho }\left( x\right) >c>-c>-\xi _{\rho }\left(
x\right) $ and the Oseledets decomposition extends continuously to $X$.

\item $\rho $ is a differentiable point of $\rho \mapsto \int_{X}\xi _{\rho
}\left( x\right) \nu \left( dx\right) $.
\end{enumerate}
\end{proposicao}

This fact extends to any real rank $1$ Lie group $G$, which means that $\dim 
\mathfrak{a}=1$. For such a group there is a unique flag manifold $\mathbb{F}
$ which is diffeomorphic to a sphere. An open Bruhat cell in $\mathbb{F}$ is
the complement of a point hence any proper subset $C=\mathrm{cl}\left( 
\mathrm{int}C\right) \subset \mathbb{F}$ fullfills the conditions of
Propositions \ref{propcompression} and \ref{propcompressionro}.

\subsection{Symplectic group}

Let $G=\mathrm{Sp}(n,\mathbb{R})=\{g\in \mathrm{Sl}(2n,\mathbb{R}%
):g^{T}Jg=J\}$ be the symplectic group where $J=%
\begin{pmatrix}
0 & I_{n} \\ 
-I_{n} & 0%
\end{pmatrix}%
$. Its Lie algebra $\mathfrak{sp}(n,\mathbb{R})$ is given by matrices 
\begin{equation*}
\begin{pmatrix}
A & B \\ 
C & -A^{T}%
\end{pmatrix}%
\qquad B-B^{T}=C-C^{T}=0
\end{equation*}%
where a choice of $\mathfrak{a}$ is the subalgebra of matrices 
\begin{equation*}
\begin{pmatrix}
D & 0 \\ 
0 & -D%
\end{pmatrix}%
\end{equation*}%
with $D$ diagonal. If $\lambda _{i}$ stands for the map which assocites to $%
D $ its $i$-th diagonal entry then the roots are $\lambda _{i}-\lambda _{j}$%
, $i\neq j$, and $\lambda _{i}+\lambda _{j}$. A standar choice of simple
roots is 
\begin{equation*}
\Sigma =\{\lambda _{1}-\lambda _{2},\ldots ,\lambda _{n-1}-\lambda
_{n},2\lambda _{n}\},
\end{equation*}%
whose bassic weights are $\omega _{i}=\lambda _{1}+\cdots +\lambda _{i}$, $%
i=1,\ldots ,n$.

The flag manifolds of $\mathrm{Sp}(n,\mathbb{R})$ are formed by flags of
isotropic subespaces. We focus on the Grassmannian Lagrangean $\mathrm{LGr}%
_{n}(2n)$, which is the set of Lagrangean subspaces of $\mathbb{R}^{2n}$.
The subset $\Theta _{\mathrm{Lag}}\subset \Sigma $ corresponding to $\mathrm{%
LGr}_{n}(2n)$ is $\Theta _{\mathrm{Lag}}=\{\lambda _{1}-\lambda _{2},\ldots
,\lambda _{n-1}-\lambda _{n}\}$.

Let $Q$ be the quadratic form on $\mathbb{R}^{2n}$ with matrix $%
\begin{pmatrix}
0 & I_{n} \\ 
I_{n} & 0%
\end{pmatrix}%
$ and define the semigroup 
\begin{equation*}
S_{Q}=\{g\in \mathrm{Sp}(n,\mathbb{R}):Q(g(v))\geq Q(v)\}.
\end{equation*}

Cocycles with values in this semigroup were considred before (see Wojtkoweki 
\cite{woj} and Barreira-Pesin \cite{barrpe}) and estimates for the Lyapunov
exponents were obtained.

To apply Corollary \ref{corsemigroupres} we note that the flag type of $%
S_{Q} $ is $\mathrm{LGr}_{n}(2n)$ as follows by general results on
compression semigroups of Hilgert-Neeb \cite{hilgnee} and \cite{smmax}.
Hence if a cocycle $\rho $ takes values in 
\begin{equation*}
\mathrm{int}S_{Q}=\{g\in \mathrm{Sp}(n,\mathbb{R}):Q(g(v))>Q(v)\}
\end{equation*}%
then the $\mathfrak{a}$-Lyapunov spectrum of the ensuing flow $\phi $ has
the form 
\begin{equation*}
\Lambda _{\mathrm{Ly}}\left( \phi \right) =\left( 
\begin{array}{ll}
D\left( \phi \right) & 0 \\ 
0 & -D\left( \phi \right)%
\end{array}%
\right)
\end{equation*}%
with $D\left( \phi \right) =\mathrm{diag}\{\chi _{1}\left( \phi \right)
,\ldots ,\chi _{n}\left( \phi \right) \}$, $\chi _{1}\left( T\right) \geq
\cdots \geq \chi _{n}\left( T\right) $ such that $2\chi _{n}\left( \phi
\right) >0$, by Corollary \ref{corsemigroupres} (2). The first statement of
the same corollary ensures that $\chi _{1}\left( \phi \right) +\cdots +\chi
_{n}\left( \phi \right) $ is differentiable as a function of the cocycle.

\end{document}